%% file: ckpaper.tex
\documentclass[12pt]{article}
\usepackage{amsfonts}
\usepackage{amssymb}
\usepackage[mathscr]{eucal}
\usepackage{amsmath}
\usepackage[dvips]{epsfig,graphics}
\usepackage[all]{xy}

\input{ckpaperchapters.tex}

\begin{document}
\setcounter{page}{1}
\title{Vassiliev invariants and the cubical knot complex}

\author{Ilya Kofman\\  
{\em {\small Department of Mathematics, University of Maryland }}\\ {\em {\small College Park,
MD 20742}}\\ \texttt{\small ikofman@math.umd.edu} \\ \\ Xiao-Song Lin\\ {\em {\small
Department of Mathematics, University of California, Riverside }}\\ {\em {\small Riverside, CA
92521}}\\ \texttt{\small xl@math.ucr.edu}}

\date{September 27, 2000}

\maketitle

\begin{abstract}
We construct a cubical CW-complex $CK(M^3)$ whose rational cohomology algebra contains
Vassiliev invariants of knots in the 3-manifold $M^3$.  We construct
$\overline{CK}(\mathbf{R}^3)$ by attaching cells to $CK(\mathbf{R}^3)$ for every degenerate
1-singular and 2-singular knot, and we show that $\pi_1(\overline{CK}(\mathbf{R}^3))=1$ and
$\pi_2(\overline{CK}(\mathbf{R}^3))=\mathbf{Z}$. We give conditions for Vassiliev invariants
to be nontrivial in cohomology.  In particular, for $\mathbf{R}^3$ we show that $v_2$ uniquely
generates $H^2(CK,D)$, where $D$ is the subcomplex of degenerate singular knots.  More
generally, we show that any Vassiliev invariant coming from the Conway polynomial is
nontrivial in cohomology.  The cup product in $H^*(CK)$ provides a new graded commutative
algebra of Vassiliev invariants evaluated on ordered singular knots. We show how the cup
product arises naturally from a cocommutative differential graded Hopf algebra of ordered
chord diagrams.
\end{abstract}

\tableofcontents

\section{Introduction}

In \c{V}, Vassiliev constructed a new series of numerical invariants of knots, now called
\emph{Vassiliev invariants}.  Shortly afterwards, it was shown that Vassiliev invariants can
be described combinatorially by axioms and initial data, and therefore the coefficients of
certain power-series expansions of all \emph{quantum invariants}, such as Jones-type polynomials,
are Vassiliev invariants (see \c{BL,BN,B}).  Unlike quantum
invariants, Vassiliev invariants are constructed by methods from classical topology.  Let
$\mm$ be the space of all smooth maps from $S^1$ into $\S$, and let the \emph{discriminant}
$\Si$ of $\mm$ be the set of maps which are not embeddings.  The components of $\mm \setminus
\Si$ are in one-one correspondence with knot types, so a numerical knot invariant is a class
in $H^0(\mm \setminus \Si;\Q)$. Even though $\mm$ is infinite-dimensional, Vassiliev was able to apply
Alexander duality by constructing a system of finite-dimensional approximations $\Gamma^d \subset \mm$,
along with auxiliary spaces $\s$ whose homology coincides with that of the discriminant $\Gamma^d \cap \Si$.

In the past decade, the focus of interest in Vassiliev invariants has been on their rich
algebraic structures (see \c{BN}).  Unfortunately, Vassiliev's original formulation has not
played a significant role in efforts to understand these structures.  On the other hand,
although Vassiliev invariants can be completely described by combinatorial data,
the methods of combinatorial and algebraic topology have failed to yield a proof of
Kontsevich's ``fundamental theorem of Vassiliev invariants'' \c{Kont, BNS}.

In this paper, we construct a cubical CW-complex $CK$, resembling Vassiliev's space $\s$,
whose rational cohomology algebra contains Vassiliev invariants.  Starting with an algebraic
chain complex of knots with ordered double points or base point in any oriented 3-manifold
$\M$, we geometrically realize the chain complex as the CW-complex $CK(\M)$.  For $\R^3$, we
show that any Vassiliev invariant coming from the Conway polynomial is nontrivial in
cohomology so, for example, the powers $v_2^n$ are nontrivial generators of
$H^{2n}(CK(\R^3);\Q)$.  In particular, we show that $v_2$ uniquely generates
$H^2(CK(\mathbf{R}^3),D)$, where $D$ is the subcomplex of degenerate singular knots.

The cup product in $H^*(CK)$ provides a new graded commutative algebra of Vassiliev invariants
evaluated on ordered singular knots.  From the algebraic point of view, we show how the cup
product arises naturally from a bialgebra structure on \emph{ordered} chord diagrams.  The
resulting cocommutative differential graded Hopf algebra (DGH) also has an interesting quotient
DGH of \emph{oriented} chord diagrams, which is commutative and cocommutative in the
graded sense, but it is open whether this quotient is nontrivial.

On the topological side, by a deep theorem of Quillen \c{Q}, the algebra of Vassiliev invariants under the cup product is the rational cohomology ring of a simply-connected pointed topological space.  We construct a knot complex $\tck(\R^3)$ which seems very close to such a space for even Vassiliev invariants.  We attach cells to $CK(\R^3)$ for every \emph{degenerate} 1-singular and 2-singular knot, and show that $\pi_1(\tck(\R^3))=1$ and $\pi_2(\tck(\R^3))=\Z$.  The proof is an
application of the conceptual extension of Vassiliev invariants to more general 3-manifolds by
studying Map$(S^1 \times D^2,\M)$ in ``almost general position'' \c{L}. Kalfagianni \c{Kalf}
further extended these techniques to show that Vassiliev invariants exist for knots in a large
class of 3-manifolds.


This paper is organized as follows.  In Section \ref{sec2}, we construct a commutative diagram
of chain complexes of singular knots with additional structures, such as ordered double
points, base point, and orientation coming from the ordering (Theorem \ref{diagthm}). For two
of these chain complexes, we construct their geometric realizations, the knot complexes $CK_b$
and $CK_0$ (Theorem \ref{CW}).  Let $CK$ denote either complex.  In Section \ref{homotopy}, we
construct $\tck$ and show that $\pi_1(\tck(\R^3))=1$ (Corollary \ref{pi1}) and
$\pi_2(\tck(\R^3))=\Z$ (Theorem \ref{pi2}).  In Section \ref{mark}, we show that a subcomplex
of markings in the knot complex $CK_b(\S)$, resembling a discriminant, can be associated with
the cellular decomposition of $\s$ described in \c{V}. In Section \ref{sec5}, we discuss the
relationship between Vassiliev invariants and the cohomology of $CK$.  Our main result is that
if a Vassiliev invariant is nontrivial in the cohomology of $CK$, then it is of even order,
and the cohomology dimension must equal the order (Corollary \ref{Cor2k}).  In Subsection
\ref{mirror}, we define a topological map $F:CK(\R^3)\to CK(\R^3)$ obtained from the mirror
map. In Subsection \ref{cup}, we discuss the cup product in $H^*(CK)$. In Section \ref{Hopf},
we define Hopf algebra structures for ordered and oriented chord diagrams, and show how the
cup product in $H^*(CK_0)$ for Vassiliev invariants arises from a Hopf algebra of ordered
chord diagrams (Theorem \ref{Wcup}).

\subsection*{\sc Acknowledgments}
This work includes part of the first author's Ph.D. dissertation at the University of Maryland, College Park, under the direction of William Goldman and Yongwu Rong (GWU).  The first author would also like to thank Dror Bar-Natan, J\'ozef Przytycki, James Schafer, Adam Sikora, and Akira Yasuhara for helpful discussions, and acknowledge support by the NSF and by the UMCP Mathematics Dissertation Fellowship.  The second author acknowledges support by the NSF.  This project first took shape in the summer of 1998 at the University of Iowa, where the second author gave a course in the IMA graduate summer school ``Topology of Manifolds''.  We would like to thank the University of Iowa mathematics department for providing a stimulating environment during that period.

\input{ckpapertext.tex}

\addcontentsline{toc}{section}{Bibliography}
\bibliography{research}
\bibliographystyle{plain}

\end{document}

%% file: ckpaperchapters.tex
\newcommand{\Z}{\mathbf{Z}}
\newcommand{\R}{\mathbf{R}}
\newcommand{\Q}{\ensuremath{\mathbf{Q}}}

\newcommand{\V}{\mathbb{V}}
\newcommand{\eop}{\hfill \ensuremath{\square} \\}
\newcommand{\pf}{\noindent {\it Proof:} \quad}
\renewcommand{\c}[1]{\cite{#1}}

\newcommand{\f}[1]{\frac{#1}{2}}

\newcommand{\s}{\ensuremath{\sigma}}
\renewcommand{\t}{\ensuremath{\tau}}
\newcommand{\Si}{\ensuremath{\Sigma}}
\newcommand{\tphi}{\ensuremath{\tilde{\Phi}}}
\renewcommand{\ss}{\ensuremath{\mathcal{\it{\zeta}}_\Phi}}
\newcommand{\tss}{\ensuremath{\mathcal{\it{\zeta}}_{\tilde{\Phi}}}}

\newcommand{\M}{\ensuremath{M^{3}}}
\newcommand{\mm}{\ensuremath{\mathscr{M}}}
\renewcommand{\S}{\ensuremath{S^{3}}}

\newcommand{\Kx}{\ensuremath{K_{\times}}}
\newcommand{\Kp}{\ensuremath{K_{+}}}
\newcommand{\Km}{\ensuremath{K_{-}}}
\newcommand{\Kn}{\ensuremath{K_{1 \ldots n}}}
\newcommand{\KM}{\ensuremath{K_{1 \ldots m}}}
\newcommand{\Ksn}{\ensuremath{K_{\sigma(1 \ldots n)}}}
\newcommand{\Ktn}{\ensuremath{K_{\tau(1 \ldots n)}}}
\newcommand{\Kwn}{\ensuremath{K_{1 \ldots n}^\omega}}

\newcommand{\Kdots}{\ensuremath{K_{\times \ldots \times}}}
\newcommand{\xdots}{{\times \ldots \times}}
\newcommand{\di}{\ensuremath{\mathrm{d}_{i}}}
\renewcommand{\d}[1]{\ensuremath{\mathrm{d}_{#1}}}
\newcommand{\del}{\ensuremath{\partial}}
\newcommand{\lr}{\ensuremath{\, \overset{1-1}{\longleftrightarrow} \,}}
\newcommand{\n}{\ensuremath{n}-}
\newcommand{\Ra}{\ensuremath{\Longrightarrow}}
\renewcommand{\f}{\ensuremath{f_{\Kn}}}
\newcommand{\g}{\ensuremath{g_{K_{1 \ldots n}}}}

\newcommand{\e}{\ensuremath{\varepsilon}}
\renewcommand{\r}{\ensuremath{\varrho}}
\newcommand{\tck}{\ensuremath{\overline{CK}}}  
\renewcommand{\th}{\ensuremath{\tilde{H}}}
\newcommand{\Dn}{\ensuremath{D_{1 \ldots n}}}
\newcommand{\Dsn}{\ensuremath{D_{\sigma(1 \ldots n)}}}

\renewcommand{\aa}{\ensuremath{\mathcal{A}}}
\newcommand{\dd}{\ensuremath{\mathfrak{D}}}
\newcommand{\DD}{\ensuremath{\mathcal{D}}}

\newcommand{\D}{\ensuremath{\Delta}}

\newtheorem{thm}{Theorem}[section] 
\newtheorem{defn}[thm]{Definition}
\newtheorem{lem}[thm]{Lemma}
\newtheorem{cor}[thm]{Corollary}
\newtheorem{ex}[thm]{Example}
\newtheorem{prop}[thm]{Proposition}

\newtheorem{rmk}[thm]{Remark}
\newtheorem{conj}[thm]{Conjecture}
\newtheorem{ques}[thm]{Question}


%% file: ckpapertext.tex
\section{Constructing the knot complexes $CK(\M)$} \label{sec2}

Let $\M$ be an oriented 3-manifold.  Let $X_{n}$ be the set of equivalence classes of oriented
knots with $n$ double points (or ``\n singular knots'') in $\M$.  These are immersions $S^{1}
\to \M$ which are embeddings except for $n$ transverse self-intersections, equivalent up to
rigid vertex isotopy.  We will denote elements of $X_{n}$ by $K_{\underbrace{\xdots}_{n}}$.
When we consider ordered double points, we will denote an ordering up to even permutations as an
``orientation,'' but this is not the orientation of $S^{1} \subset \M$.  To distinguish our terms, we
will use \emph{direction} to refer to the usual orientation of $S^1\subset\M$.

\begin{figure}
\begin{center}
\scalebox{.45}{\includegraphics{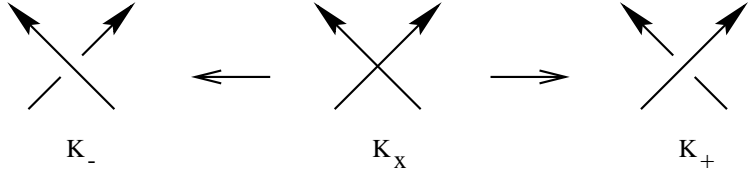}}
\caption{Resolving a double point} \label{skeinfig}
\label{Kx.eps}
\end{center}
\end{figure}

The choice of direction of a singular knot is unimportant for resolving its double points
because we can resolve any double point in two canonical ways by considering the orientation
of the ambient manifold $\M$.  See Figure \ref{skeinfig}.  (This is not true for links, and a
similar construction for links would not be invariant under changes of direction.) We call a
singular knot \emph{degenerate} if for some singular crossings, their positive and negative
resolutions determine the same knot type.

\begin{rmk} \rm
A singular crossing in a knot diagram is called \emph{nugatory} if there exists a separating
$S^2$ which intersects the diagram in only one point.  Clearly, such a crossing is degenerate.
It seems to be unknown whether there exist degenerate crossings which are not nugatory.
\end{rmk}

\subsection{Additional structures on singular knots} \label{structures}

\paragraph{\n singular knots with ordered double points}
Let $X_n^0$ be the set of equivalence classes of elements $\Kdots \in X_n$ with double points ordered from $1$ to $n$.  Let $X_0^0=X_0$.  We will denote elements of $X_n^0$ by $\Kn$.  For any $\Kn\in X_n^0$ and any $\sigma\in S_n$, let $\Ksn$ denote the knot obtained by permuting the ordering of the double points of $\Kn$.  In $X_n^0$, $\Kn=\Ksn$ if there exists a direction-preserving rigid-vertex isotopy $ \theta:\M\to\M$, such that $\theta(\Kn)= \Ksn$ by mapping each double point to the one with the same label.  In this case, we say $\theta$ is \emph{order-preserving}.  Let $C_n(X^0)$ denote the free abelian group with coefficients in $\Q$ generated by elements $\Kn\in X_n^0$.

\paragraph{\n singular knots with a base point}
Let $X^b_n$ be the set of equivalence classes of elements $\Kdots \in X_n$ with a base point $b$ which is apart from all the double points.  Here, two based knots $(K,b),(K',b')$ are equivalent if there exists a direction-preserving rigid-vertex isotopy $\theta:\M\to\M$ such that $(\theta(K),\theta(b))=(K',b')$, and the base point remains apart from the double points during the isotopy.  In this case, we say $\theta$ \emph{preserves base points}.  For any knot in $X^b_{n}$, we can order its double points from $1$ to $n$ around the knot starting from the base point.  (This requires the knot to be oriented.)  We will denote based \n singular knots with this natural ordering by $(\Kn,b)$.  Let $C_n(X^b)$ denote the free abelian group with coefficients in $\Q$ generated by elements $(\Kn,b)\in X_n^b$.

\paragraph{\n singular knots with ordering and base point}
Let $X_n^{b,0}$ be the set of elements in $X^b_n$ with the ordering of the double points permuted by any $\sigma\in S_n$.  Let $X_0^{b,0}=X^b_0$.  As sets $X_n^{b,0} \cong X^b_n \times S_n$.  We will denote elements of $X_n^{b,0}$ by $(\Ksn,b)$.  Let $C_n(X^{b,0})$ denote the free abelian group with coefficients in $\Q$ generated by elements $(\Ksn,b) \in X_n^{b,0}$. \\

From the definitions above, $X^0_n=X_n^{b,0}/\sim_0$, where $(\Ksn,b) \sim_0 (\Ktn,b')$ if there exists a direction-preserving rigid-vertex isotopy which maps each double point to the one with the same label, but need not preserve base points.  Let $p_0: C_n(X^{b,0})\to C_n(X^0)$ be the quotient map by $\sim_0$.
Let $i: C_n(X^b) \to C_n(X^0)$ be the map $i(\Kn,b)=p_0(\Kn,b)$.

Because $X_n^{b,0} \cong X^b_n \times S_n$, it follows that $C_n(X^b)\cong C_n(X^{b,0})/\sim_b$, where
\[ (\Kn,b) \sim_b \text{\rm sign} (\sigma) (\Ksn,b) \; \text{ for any } \sigma \in S_n \]
Let $p_b: C_n(X^{b,0})\to C_n(X^b)$ be the map $p_b(\Ksn,b)=\text{\rm sign} (\sigma)(\Kn,b)$.

Let $C_n(X^\omega)$ be the abelian group with coefficients in $\Q$ generated by elements $\Kn \in X_n^0$, subject to the following relations:
\[ \Kn \sim_\omega \text{\rm sign} (\sigma) \Ksn \; \text{ for any } \sigma \in S_n \]
Let $p_\omega: C_n(X^0)\to C_n(X^\omega)$ be the quotient map by $\sim_\omega$.

\begin{prop}  The following diagram commutes:
\begin{equation} \label{CnX}
\xymatrix{ {} & \ar[dl]_{p_b} {C_n(X^{b,0})} \ar[dr]^{p_0} & {}  \\
{C_n(X^b)} \ar[dr]_{p_\omega\circ i}  & {} & \ar[dl]^{p_\omega} {C_n(X^0)} \\
{} & {C_n(X^\omega)} & {} }
\end{equation}
\end{prop}
\pf
\[ p_\omega\circ i\circ p_b(\Ksn,b)=p_\omega(\text{\rm sign} (\sigma)\Kn)=p_\omega(\Ksn)=p_\omega\circ p_0(\Ksn,b) \]
\eop

\begin{defn} \label{orientable}
For $n \geq 2, \; \Kdots \in X_n$\, is \emph{non-orientable} if there is an odd permutation $\sigma\in S_n$, such that $\Kn=\Ksn$ for some ordering of the double points.  Otherwise, we say $\Kdots$ is \emph{orientable}.  Let $\Omega_n \subset X_n,\;\Omega_n^0 \subset X_n^0,$ and $\Omega_n^b \subset X_n^b$ be the respective subsets of orientable $n$-singular knots.
\end{defn}

Since $X_n$ equals $X_n^0$ modulo the action of $S_n$ on orderings, the projection $(X_n^0/A_n)\to X_n$ is either $2-1$ or $1-1$.  Non-orientable knots are exactly those with one preimage, so they are exactly the knots killed by $\sim_\omega$, and hence are in the kernel of $p_\omega$.  If $n\in\{0,1\},\; p_\omega$ is an isomorphism.

\paragraph{\n singular knots with orientation}
For every $\Kdots\in \Omega_n$, there are two generators in $C_n(X^\omega)$, one in each coset of $\Omega_n^0/A_n$.  Let $\Kwn$ be either generator.  An \emph{orientation} $\omega$ is the union for $n \geq 2$ of choices $X_n^\omega$, where $X_n^\omega = \bigcup\Kwn$.  If $\Kdots\in \Omega_n$, we will refer to the preimage of $\Kwn$ in $X_n^{b,0}$ as knots in its orientation class, and $\omega$ will also denote the preimage of $X_n^\omega$ in $X_n^{b,0}$.

The choice of cosets $\Omega_n^0/A_n$ is unimportant because there is a chain complex isomorphism $\Gamma: C_n(X^\omega) \to C_n(X^{\omega'})$ by $\Gamma(\Kwn)=\pm \Kn^{\omega'}$ for every $\Kdots\in \Omega_n$, where the sign is positive if and only if $\omega$ and $\omega'$ impose the same ordering on $\Kdots$.

\begin{figure}
\begin{center}
\scalebox{.55}{\includegraphics{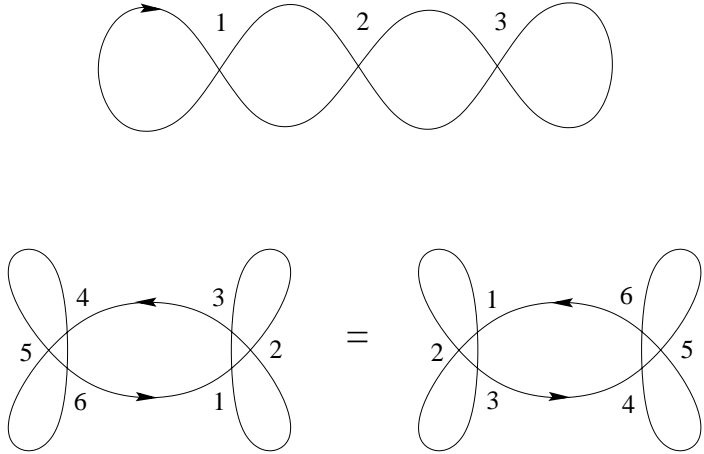}}
\caption{Non-orientable knots}
\label{nonorientable.eps}
\end{center}
\end{figure}

\begin{ex} \rm
We refer to Figure \ref{nonorientable.eps}.  If $K_{1 \ldots 6}$ is the $6$-singular knot shown, then $K_{1 \ldots 6}=K_{\sigma(1 \ldots 6)}$, where $\sigma(123456)=(456123)$ is an odd permutation.  If $\Kn$ is the degenerate knot shown (in the figure $n=3$), then $\Kn= \Ksn$, where $\sigma(1\ldots n)=(n\ldots 1)$.  Since $\sigma$ is the product of $\binom{n}{2}=n(n-1)/2$ transpositions, $\sigma$ is odd if and only if $n\equiv 2\mod 4$ or $n\equiv 3\mod 4$.  This provides an infinite family of non-orientable knots.
\end{ex}

\subsection{Constructing chain complexes of singular knots}

\noindent
Let $C_n(X)$ denote any of the groups from the previous section.  For the corresponding elements, we will suppress the notation for base points whenever possible.

\begin{defn}
Given $\Ksn \in X_n^{b,0}$, define $\di: C_n(X^{b,0}) \to C_{n-1}(X^{b,0})$ by
\[ \di \Ksn = K_{\times \ldots \underset{i}{+} \ldots \times} - K_{\times \ldots \underset{i}{-} \ldots \times} \]
such that each resolution has its ordering induced from $\Ksn$.
\end{defn}

\begin{defn} \label{deg}
For $n \geq 1, \,\, \Kdots \in X_n$\, is \emph{degenerate} if given any ordering of its double points, $\exists i, 1\leq i \leq n,$ such that $ \di \Kdots = 0.$  Let $D_n \subset X_n$ be the set of degenerate $n$-singular knots.  Let $C_n(D)$ denote the respective subgroup of $C_n(X)$ generated by degenerate \n singular knots.
\end{defn}

\begin{defn} \label{del}
Define the boundary operator $\del: C_n(X^{b,0}) \to C_{n-1}(X^{b,0})$ by
\[ \del \Ksn = \sum_{i=1}^{n} (-1)^{i+1} \di \Ksn \]
\end{defn}

Taken together, the Propositions \ref{del1}, \ref{del2}, \ref{del3}, and \ref{del4} below prove the following theorem:

\begin{thm} \label{diagthm}
With the boundary operator defined above, (\ref{CnX}) is a commutative diagram of chain complexes and chain maps.  Moreover, with the same boundary operator, we can form the relative chain complexes $C_n(X,D)=C_n(X)/C_n(D)$.
\end{thm}

\begin{prop} \label{del1}
If $(\Ksn,b) \sim_0 (\Ktn,b')$ then $\del\Ksn=\del\Ktn$, so the boundary operator is well-defined on $C_n(X^0)$.
\end{prop}
\pf \qquad Let $K_{\times\ldots\e_i\ldots\times}^\sigma$ denote the $i^{\text{th}}$ resolution of $K_{\sigma(1 \cdots n)}$ with its ordering induced from the ordering $K_{\sigma(1 \cdots n)}$.
By assumption, there is an order-preserving isotopy $\theta: \Ksn \to \Ktn$, so it induces an isotopy on resolutions: $K_{\times\ldots\e_i\ldots\times}^\sigma=K_{\times\ldots\e_i\ldots\times}^\tau$.  Therefore,
\[ \begin{array}{l}
\del\Ksn=\sum\limits_{i=1}^{n}\sum\limits_{\e\in\pm 1}(-1)^{i+1}\e_i (K_{\times\ldots\e_i\ldots\times}^\sigma) = \\ \\
\sum\limits_{i=1}^{n}\sum\limits_{\e\in\pm 1}(-1)^{i+1}\e_i (K_{\times\ldots\e_i\ldots\times}^\tau)=\del\Ktn
\end{array}  \]
\eop
\begin{prop}{$\del^2 = 0.$} \label{del2}
\end{prop}
\pf \qquad For any knot in $X_n^{b,0},\, X_n^0,\, X_n^b$, if $1 \leq i< j \leq n, \, \di\d{j} = \d{j-1}\di \,$.  Therefore,
\[ \del^2 = \sum_{i=1}^{n}(-1)^{i+1}\left(\sum_{j<i}(-1)^{j+1}\d{j}\di + \sum_{j\geq i}(-1)^{j+1}\di\d{j+1}\right) = 0 \]
since if we replace $j+1$ by $j$ in the second sum, every term appears twice with opposite sign.  Degenerate knots will result in additional zero terms.
\eop
\begin{prop}{For $n \geq 2$, if $\forall \sigma \in S_n, \; \Ksn = \text{sign} (\sigma)\Kn$, then $\del\Ksn=\text{sign} (\sigma)\del\Kn$.} \label{del3}
\end{prop}
\pf\qquad It suffices to consider the transposition $\sigma=(p,p+1) \in S_n$.
\[ \begin{array}{l}
\del\Ksn = \sum_{i=1}^{n} (-1)^{i+1} \di \Ksn \\ \\
= \sum\limits_{i\neq p,p+1}(-1)^{i+1}\di\Ksn + (-1)^{p+1}\d{p}\Ksn + (-1)^{p+2}\d{p+1}\Ksn \\ \\
= \sum\limits_{i\neq p,p+1}(-1)^{i+1}(-\di\Kn) + (-1)^{p+1}\d{p+1}\Kn + (-1)^{p+2}\d{p}\Kn \\ \\
= -\sum_{i=1}^{n} (-1)^{i+1} \di \Kn = -\del\Kn \hfill \square \\
\end{array} \]
\begin{prop} \label{del4}
Let $C_n(D)$ denote the abelian group with coefficients in $\Q$ generated by degenerate \n singular knots with ordered double points.  For $n \geq 1, \; \del(D_n) \subset C_{n-1}(D)$, where we let $C_0(D)=0$.
\end{prop}

\begin{rmk} \rm
A similar chain complex has been used to generalize skein modules \c{BFK}, and an algebraic chain complex was constructed by Hutchings \c{Hu} to prove the analogue of Kontsevich's theorem for singular braids.  However, it does not seem possible to realize these algebraic chain complexes as CW-complexes.
\end{rmk}

\subsection{Geometric realization of the chain complexes}

In this section, we construct CW-complexes $CK_b$ and $CK_0$ whose cellular chain complexes are isomorphic to the chain complexes $C_*(X^b)$ and $C_*(X^0)$, respectively.  We construct $CK_b$ and $CK_0$ by attaching \n cubes in one-one correspondence with \n singular knots, respectively with a base point or ordering:

\begin{itemize}
\item vertices ($0$-cubes) \lr oriented knot types $K \in X_0$
\item for $n \geq 1,$ \n cubes \lr \n singular knots $\Kn\in X_n^b$ or $X_n^0,$ respectively
\end{itemize}

\paragraph{Notation} Henceforth, $CK$ will denote either $CK_0$ or $CK_b$, and all unspecified statements are valid for both complexes. \\

\begin{figure}
\begin{center}
\scalebox{.45}{\includegraphics{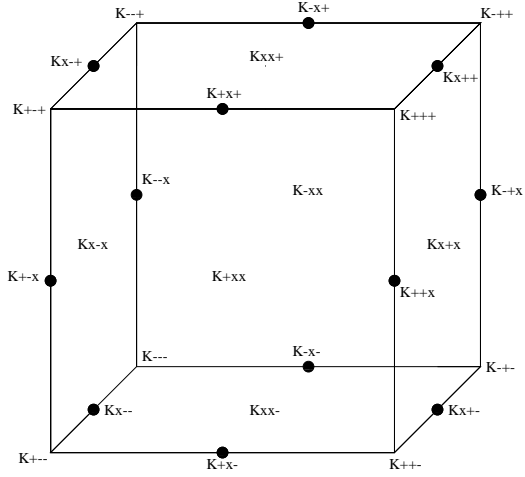}}
\caption{Labelings on the 3-cube $K_{123}$}
\label{Kxxx}
\end{center}
\end{figure}

Let $CK^{(n)}$ denote the \n skeleton of $CK$.  $CK^{(1)}$ can be constructed immediately from the formula $\del \Kx = \Kp - \Km$.  Let $I=[-1,1]$.  For any $1$-singular knot $\Kx$, the attaching map $f_{\Kx} : \del I \to CK^{(0)}$ is given by $f(-1)=\Km$ and $f(+1)=\Kp$.  The map $g_{\Kx} : I \to CK^{(1)}$ given by actually attaching the edge is called the \emph{characteristic map} of $\Kx$, and it induces an orientation of the edge from $\Km$ to $\Kp$.  Thus, any knots which differ by a single crossing change are connected by an oriented edge.  An edge joins distinct vertices if and only if the corresponding $1$-singular knot is nondegenerate.  When considered as elements of $\pi_1(\M)$, only knots in the same conjugacy class can be obtained by changing crossings, so components of $CK(\M)$ correspond bijectively with conjugacy classes of $\pi_1(\M)$.

Having constructed $CK^{(1)}$, we proceed inductively.  Let $I^n$ be the \n cube $\{\vec{x} \in \R^n : -1 \leq x_i \leq 1 \}$ with the standard orientation.  Assuming $CK^{(n-1)}$ is constructed, for every $\Kn$ we attach an \n cube by the attaching map $\f :\del I^n \to CK^{(n-1)}$.  Suppose $p+q=n$.   We label the $p$-faces of $I^n$ with $p$-singular knots obtained by resolving $q$ singularities as follows (See Figure \ref{Kxxx}):
\[ K_{\times\ldots\e_1\ldots\e_q\ldots\times} \longleftrightarrow \{(x_1,\dots,\e_1,\dots,\e_q,\dots,x_p):\e_i\in\pm 1, -1\leq x_j\leq 1 \} \]

Let $K_{1\ldots n-1}$ denote any $(n-1)$-singular resolution of $\Kn$ with its ordering induced from $\Kn$.
By assumption, for every $K_{1 \cdots n-1}$, there is a characteristic map $g_{K_{1 \cdots n-1}} : I^{n-1} \to CK^{(n-1)}$.  We define $\f :\del I^n \to CK^{(n-1)}$ by mapping each $(n-1)$-face according to its labeling.  For $i\in\{1,\dots,n\},\;\e_i\in\pm 1$, we resolve the $i^{\text{th}}$ crossing and define for $\vec{x} \in K_{\times\ldots\e_i\ldots\times} \subset \del I^n,$
\[ \f(x_1,\dots,\e_i,\dots,x_n) = g_{K_{\times\ldots\e_i\ldots\times}}(x_1,\dots,x_{i-1},x_{i+1},\dots,x_n) \]
To verify that $\f$ is well-defined on intersections, if $K_{\times\ldots\e_i\ldots\e_j\ldots\times}$ is an $(n-2)$-face of $\Kn$ then
\[ \begin{array}{l}
\f(x_1,\dots,\e_i,\dots,\e_j,\dots,x_n)
=g_{K_{\times\ldots\e_i\ldots\times}}(x_1,\dots,\hat{x_i},\dots,\e_j,\dots,x_n) \\ \\
=g_{K_{\times\ldots\e_i\ldots\e_j\ldots\times}}(x_1,\dots,\hat{x_i},\dots,\hat{x_j},\dots,x_n)
=g_{K_{\times\ldots\e_j\ldots\times}}(x_1,\dots,\e_i,\dots,\hat{x_j},\dots,x_n)
\end{array} \]
This completes the induction.

\begin{thm}{There exist CW-complexes $CK_b$ and $CK_0$ such that for $X=X^b$ and $X^0$, respectively, $C^{CW}_*(CK) \cong C_*(X)$.} \label{CW}
\end{thm}
\pf\qquad
Let $\Phi : C_*(X) \to C^{CW}_*(CK)$ by $\Phi(\Kn)=\g(I^n),\;\forall \Kn\in X_n^b$ or $X_n^0$.  We orient each \n cell by taking the orientation induced by $\g$ from the standard orientation of $I^n\subset\R^n$.  Since $\Phi$ maps generators bijectively and these groups are free abelian, $\Phi$ is an isomorphism of groups.  It remains to show $\Phi\circ\del=\del\circ\Phi$.  The boundary map for the cellular complex of a CW-complex is given as follows:  for any \n cell $\s, \, \del \s = \sum_{\t}[\s:\t]\t$, where $[\s:\t]$ is the incidence number with respect to the chosen orientations.  By construction, $[\Phi(\Kn):\t]\neq 0$ only if there exists some resolution $K_{\times\ldots\e_i\ldots\times}$ of $\Kn$ such that $\Phi(K_{\times\ldots\e_i\ldots\times})=\t$.
\[ \begin{array}{l}
\del \Phi(\Kn)= \sum_{\t}[\Phi(\Kn):\t]\t \\ \\
=\sum\limits_{i=1}^{n}\sum\limits_{\e\in\pm 1}[\Phi(\Kn):\Phi(K_{\times\ldots\e_i\ldots\times})]\Phi(K_{\times\ldots\e_i\ldots\times}) \\ \\
=\sum\limits_{i=1}^{n}\sum\limits_{\e\in\pm 1}[\g(I^n):g_{K_{\times\ldots\e_i\ldots\times}}(I^{n-1})]g_{K_{\times\ldots\e_i\ldots\times}}(I^{n-1}) \\ \\
=\sum\limits_{i=1}^{n}\sum\limits_{\e\in\pm 1}[\g(I^n):\g(\{\vec{x}\in\del I^n:x_i=\e_i\})]g_{K_{\times\ldots\e_i\ldots\times}}(I^{n-1}) \\ \\
= \sum\limits_{i=1}^{n}\sum\limits_{\e\in\pm 1}\left( (-1)^{i+1}\e_i\right) g_{K_{\times\ldots\e_i\ldots\times}}(I^{n-1}) \\ \\
=\sum\limits_{i=1}^{n}(-1)^{i+1}\left( \sum\limits_{\e\in\pm 1}\e_i g_{K_{\times\ldots\e_i\ldots\times}}(I^{n-1})\right) \\ \\
=\sum\limits_{i=1}^{n}(-1)^{i+1}\Phi(\di \Kn) = \Phi\del\Kn \\
\end{array} \]
\eop

\begin{prop} \label{f}
For $CK_0,$ if $\s \in A_n, \;  f_{\Ksn}\simeq\f : \del I^n \to CK_0^{(n-1)}$.
\end{prop}

\pf\qquad For any transposition $\s=(j,j+1) \in S_n$, the labeling for $\Ksn$ is obtained by reflecting the labeled \n cube corresponding to $\Kn$ in the plane $x_j=x_{j+1}$.  For any $\s \in S_n$, the labeling for $\Ksn$ is obtained by a product of reflections, one for each transposition of $\s$.  We denote this product of reflections by $\r_\s : I^n \to I^n$.  Let $K_{\times\ldots\e_i\ldots\times}^\s$ denote the $i^{\text{th}}$ resolution of $\Ksn$ with its ordering induced from the ordering $\Ksn$.  Then $\r_\s(K_{\times\ldots\e_i\ldots\times})=K_{\times\ldots\e_{\s(i)}\ldots\times}^\s$, where the same crossing is resolved on both sides.  Let $\r_\tau:I^{n-1}\to I^{n-1}$ be the product of reflections, one for each transposition of $\tau\in S_{n-1}$, which permutes the induced ordering of $K_{\times\ldots\e_{\s(i)}\ldots\times}^\s$ to $K_{\times\ldots\e_i\ldots\times}$.  For any $1 \leq i \leq n$, suppose $\vec{x} \in K_{\times\ldots\e_i\ldots\times} \subset \del I^n$.
\[ \begin{array}{l}
f_{\Ksn}\circ\r_\s(x_1,\dots,\e_i,\dots,x_n)=f_{\Ksn}(x_{\s(1)},\dots,\e_{\s(i)},\dots,x_{\s(n)}) \\ \\
=g_{K_{\times\ldots\e_{\s(i)}\ldots\times}^\s}(x_{\s(1)},\dots,\widehat{x_{\s(i)}},\dots,x_{\s(n)}) \\ \\
=g_{K_{\times\ldots\e_{\s(i)}\ldots\times}^\s}\circ\r_\tau(x_1,\dots,\hat{x_i},\dots,x_n) \\ \\
=g_{K_{\times\ldots\e_i\ldots\times}}(x_1,\dots,\hat{x_i},\dots,x_n) \\
\end{array} \]
Therefore, $f_{\Ksn}\circ\r_\s=\f$.  If $\s \in A_n,$ then $\r_\s$ is the product of an even number of reflections, so $f_{\Ksn}\simeq\f$.
\eop


\begin{rmk} \rm
Neither $CK_b$ nor $CK_0$ is locally finite.  For example, there are infinitely many knots with unknotting number $1$, so each one would be connected by an edge to the unknot.  Since a CW-complex is metrizable if and only if it is locally finite, the metric topology on $CK$ induced from the \emph{Gordian metric} on the space of knots \c{Mu}, where every edge has length $=1$, is strictly weaker than the CW topology.
\end{rmk}

\section{Homotopy of $\tck(\R^3)$} \label{homotopy}

We define $\tck$ to be the space obtained from $CK$ by attaching a 2-cell to every 1-cell
corresponding to a degenerate 1-singular knot, and attaching a 3-cell to every 2-sphere
corresponding to a degenerate 2-singular knot.   To be precise, for $i=1,2$, let $D_i$ denote the corresponding degenerate subset of $X_i^b$ or $X_i^0$.  Let $A=A_1\cup A_2$, where
$A_1=\cup_{\alpha_1\in D_1} D_{\alpha_1}^2$ and $A_2=\cup_{\alpha_2\in D_2} D_{\alpha_2}^3$.
Then $\tck=CK\cup_fA$, where $f_{\alpha_1} : \del D_{\alpha_1}^2 \to g_{\alpha_1}(I) \;
\forall \alpha_1 \in D_1$, and $f_{\alpha_2} : \del D_{\alpha_2}^3 \to g_{\alpha_2}(I^2) \;
\forall \alpha_2 \in D_2$. The map $f_{\alpha_2}$ will be more fully explained in subsection
\ref{pi2section}.
In this section, we will suppress the notation for ordering and base points whenever possible.

\subsection{$\pi_1(\tck(\R^3))=1$} \label{pi1section}

Although Theorem \ref{lin} is true for any 3-manifold, the proof of Corollary \ref{pi1} requires that $\pi_1(\M)=\pi_2(\M)=1$.  We consider a map $\Phi:S^1 \times D^2 \to \M$ as a family of piecewise linear maps $\{\phi_x:S^1 \to \M| \; x \in D^2\}$.

\begin{defn}
$\ss=$ closure\{$x \in D^2| \; \phi_x $ is not an embedding\}
\end{defn}
$\ss$ is a sub-polyhedron of $D^2.$  A 1-dimensional polyhedron is called \emph{neat} if it intersects $\del D^2$ in only finitely many boundary vertices.

\begin{thm}[Proposition 2.1 \c{L}]\label{lin}
A map $\Phi:S^1 \times D^2 \to \M$ can always be changed by an arbitrarily small perturbation so that $\ss$ is a neat 1-dimensional sub-polyhedron of $D^2$.  Moreover, we have:
\begin{enumerate}
\item if $x, x' \in D^2$ belong to the same component of $(D^2\setminus\ss))$ or \\
$(\ss\setminus$\{interior vertices\}$)$, then $\phi_x$ and $\phi_{x'}$ are ambient isotopic;
\item interior vertices of $\ss$ are of valence $4$ or $1$;
\item if $x \in \ss$ lies in an edge or is a boundary vertex, then $\phi_x$ has exactly one transverse double point;
\item if $x \in \ss$ is an interior vertex of valence $4$, then $\phi_x$ has exactly two transverse double points;
\item if $x \in \ss$ is an interior vertex of valence $1$, then $\phi_x$ is an embedding ambient isotopic to nearby embeddings.
\end{enumerate}
The resulting map $\Phi:S^1 \times D^2 \to \M$ is said to be in \emph{almost general position}.  Moreover, if $\Phi|S^1 \times \del D^2$ is already in almost general position, then the perturbation can be made relative to $\del D^2$.
\end{thm}

\begin{figure}
\begin{center}
\scalebox{.45}{\includegraphics{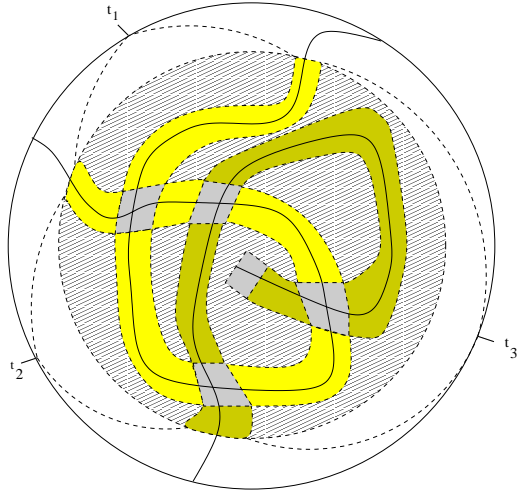}}
\caption{$\eta_\Phi$ subdivided into rectangles, and regions for extending $\Psi$ by homotopies $F_\alpha$}
\label{pi1.eps}
\end{center}
\end{figure}

\begin{cor}{If $\pi_1(\M)=\pi_2(\M)=1$, then $\pi_1(\tck(\M)) = 1.$}\label{pi1}
\end{cor}


\pf \qquad
Any loop $\gamma$ in $\tck$ is homotopic to a loop in the $1$-skeleton.  We will define a map $\Psi: D^2 \to \tck$, such that $\Psi|\del D^2 = \gamma$.  We can assume that $\gamma(t_i) \subset \tck^{(0)}$ for only finitely many $t_i \in S^1$, and otherwise $\gamma(\theta) \subset \tck^{(1)}$.  Let $h_i:\{\theta\in S^1: t_i\leq\theta\leq t_{i+1}\}\to I$ by extending linearly from $h_i(t_i)=-1, h_i(t_{i+1})=1$, where we let $t_{N+1}=t_1+2\pi$ for some $N$.  We can define $\Phi: S^1 \times S^1 \to \M,$ such that for $t_i\leq\theta\leq t_{i+1}$,
\[ \gamma(\theta)=g_{\Phi(S^1,\theta_i)}(h_i(\theta)) \]
$\Phi(S^1,\theta)$ is a $1$-singular knot for finitely many $\theta_i \in S^1$, and is a nonsingular embedding for all other $\theta \in S^1$.  $\Phi$ is defined by first identifying each vertex $\gamma(t_i)$ with any particular embedding of that knot, $\Phi(S^1,t_i)$.  On each arc for which $\gamma(\theta)$ cannot be deformed by a homotopy to a constant map, we choose an isotopy so that the only singularity is the crossing change, which occurs at $t_i<\theta<t_{i+1}$.  Therefore, $\Phi|S^1 \times S^1$ is in almost general position.  $g_{\Phi(S^1,\theta)}$ is the characteristic map for the corresponding cell in $\tck$.

Because $\pi_1(\M)=\pi_2(\M)=1$, we can extend $\Phi$ by a homotopy to $S^1 \times D^2$ rel $S^1 \times \del D^2$.  For example in $\R^3$, we may choose chords parallel to some fixed direction, and deform the knot $\Phi(S^1,\theta_0)$ to the knot $\Phi(S^1,\theta_1)$ along the chord $\overline{\theta_0 \theta_1}$.  By Theorem \ref{lin}, we can perturb $\Phi$ to be in almost general position relative to $\del D^2$.  As $\Phi(\cdot,S^1)$ extends to $\Phi(\cdot,D^2)$, we consider $\ss$ as a graph in $D^2$.  Recall that edges of the graph parametrize immersions which are $1$-singular, except that each $4$-valent vertex corresponds to a $2$-singular immersion, and each $1$-valent vertex corresponds to a cusp.  For $CK_b,\,\Phi$ can be perturbed both to be in almost general position and to preserve base points.

Let $V$ be the set of interior $1$-valent and $4$-valent vertices of $\ss$.  We can find a closed simply connected subset $B\subset int(D^2)$ such that $V \subset B$, no components of $\ss$ are entirely in $D^2 \setminus B$, and $\forall x \in \del B\cap\ss$, an arc of $\ss \setminus B$ connects $x$ with some $\theta_i\in\del D^2$.  Let $\eta_\Phi\subset B$ be the union of thin bands which contain $\ss$ and are isotopic to a tubular neighborhood of $\ss\cap B$ in $B$.  We subdivide $\eta_\Phi$ into rectangles as follows: let $\{R_\nu: \nu\in V\}$ be the intersection or self-intersection of bands at every 4-valent vertex, together with a disjoint rectangle at every 1-valent vertex.  Then $\eta_\Phi=\{R_\nu: \nu\in V\}\cup\{R_\t: \t$ is a component of $\ss\setminus V\}$.  Let $h_\t: R_\t \to I$ be extended continuously along $\ss^\perp$, starting as follows:  suppose $\Phi(S^1,\t)=\Kx$.  If $x\in R_\t\cap\del\eta_\Phi\cap int(B)$ such that $\Phi(S^1,x)=\Km$ then $h_\t(x)=-1$; similarly $h_\t(\Kp)=+1$.  For any $\nu\in V$, let $h_\nu: R_\nu \to I^2$ be defined by $h_\nu(\del R_\t)=h_\t(\del R_\t)$ whenever $\del R_\nu\cap\del R_\t\neq\emptyset$, and then extended continuously to $R_\nu$.

We now define $\Psi: B\to\tck$ as follows: if $x \in D^2 \setminus \eta_\Phi, \; \Psi(x)=g_{\Phi(S^1,x)}(0)$; if $x \in R_\alpha, \; \Psi(x)=g_{\Phi(S^1,\alpha)}(h_\alpha(x))$.  We extend $\Psi$ to $D^2$ by homotopies $F_\alpha$ for every $R_\alpha \cap \del B$.  Suppose the point $R_\alpha \cap \del B \cap \ss$ is connected by an arc to $\theta_i\in\del D^2$.  Since we can assume $h_i$ is linear, let $F_\alpha$ be the homotopy from $h_\alpha|_{R_\alpha\cap \del B}$ to $h_i$, extending $\Psi$ to $D^2\setminus B$.  On the remaining regions of $D^2\setminus B$ which are entirely nonsingular, let $\Psi(x)=g_{\Phi(S^1,x)}(0)$.  (See Figure \ref{pi1.eps}).
\eop

\subsection{$\pi_2(\tck(\R^3))=\Z$} \label{pi2section}

\begin{figure}
\begin{center}
\scalebox{.15}{\includegraphics{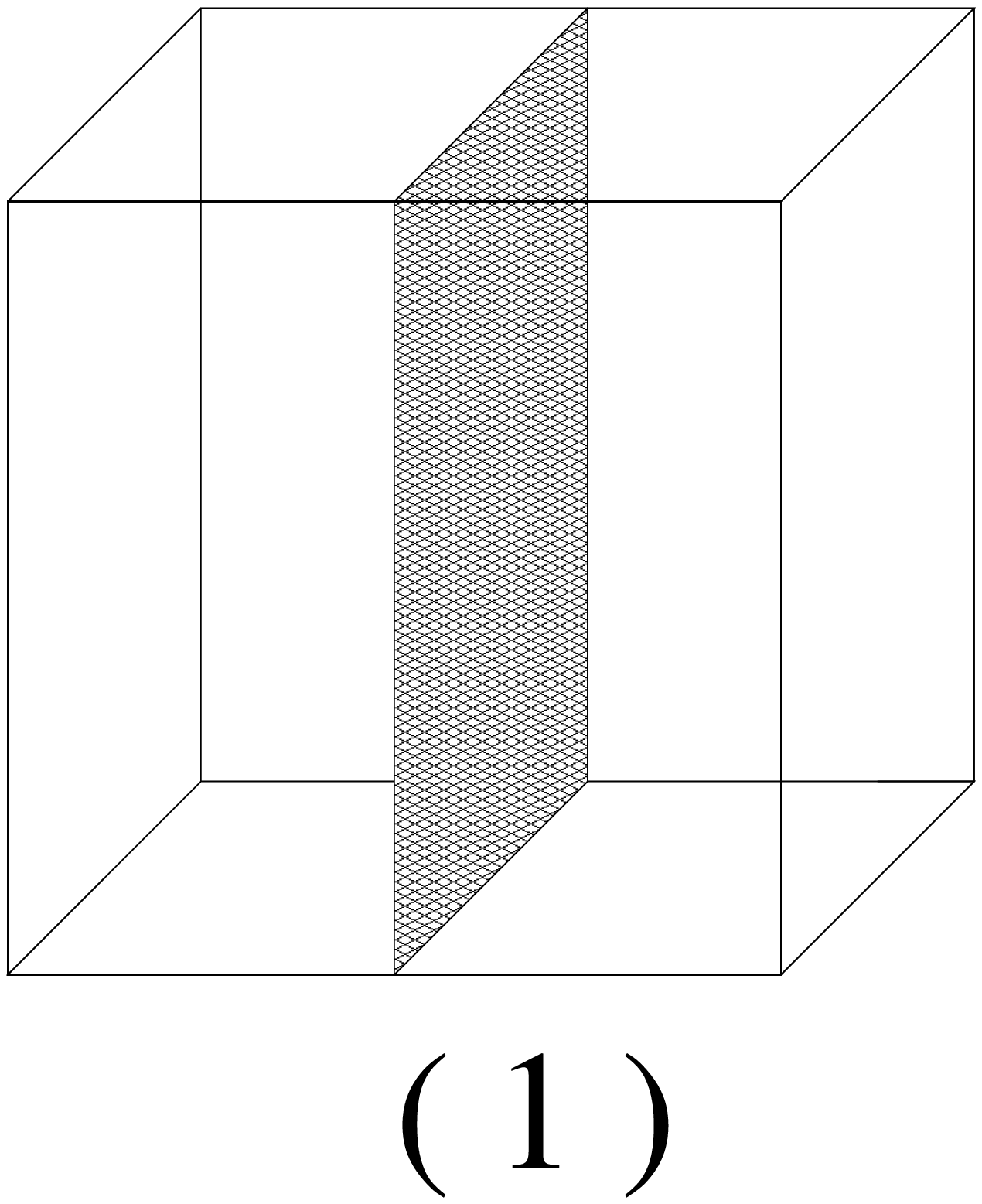}}
\scalebox{.15}{\includegraphics{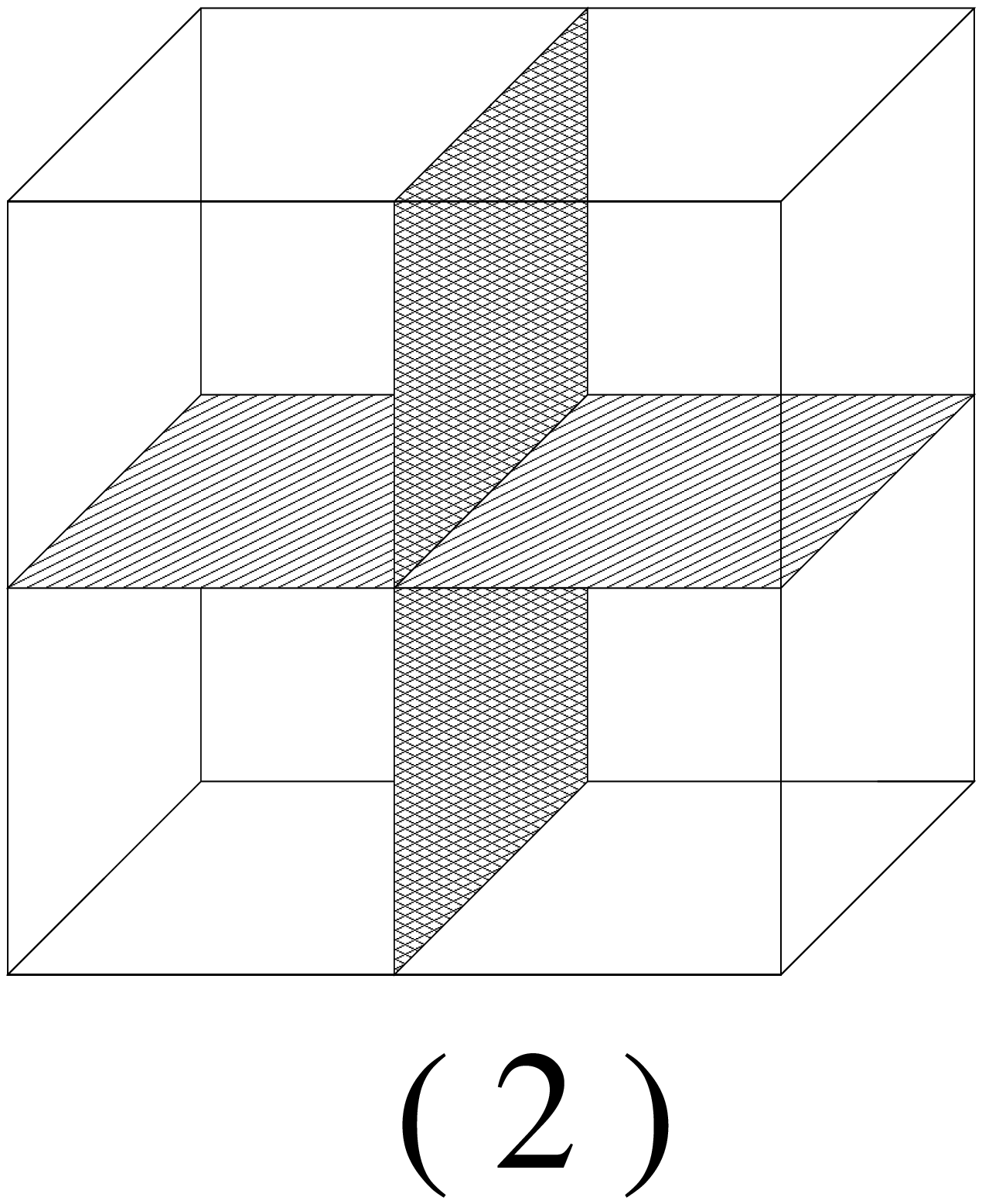}}
\scalebox{.15}{\includegraphics{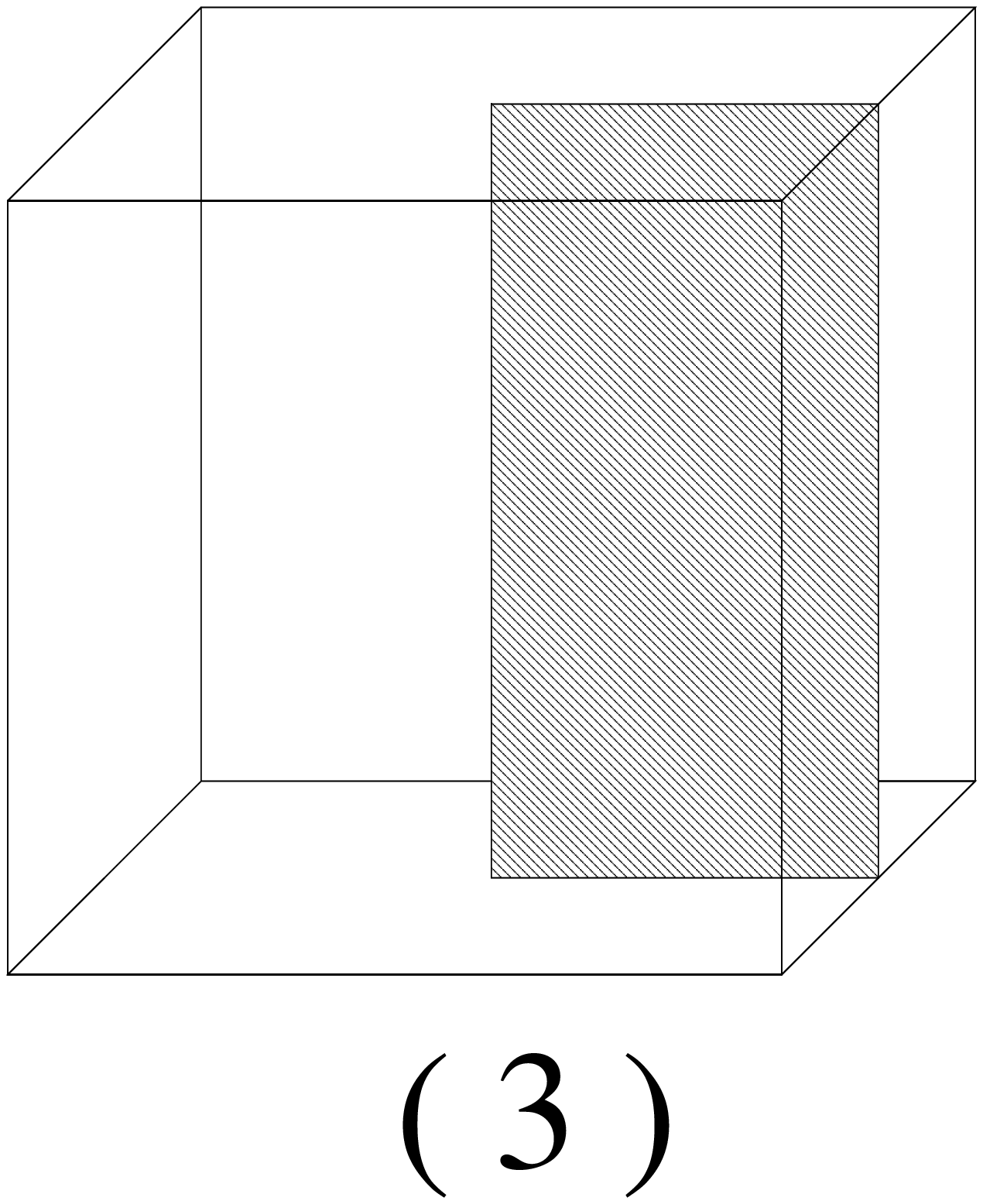}}
\scalebox{.15}{\includegraphics{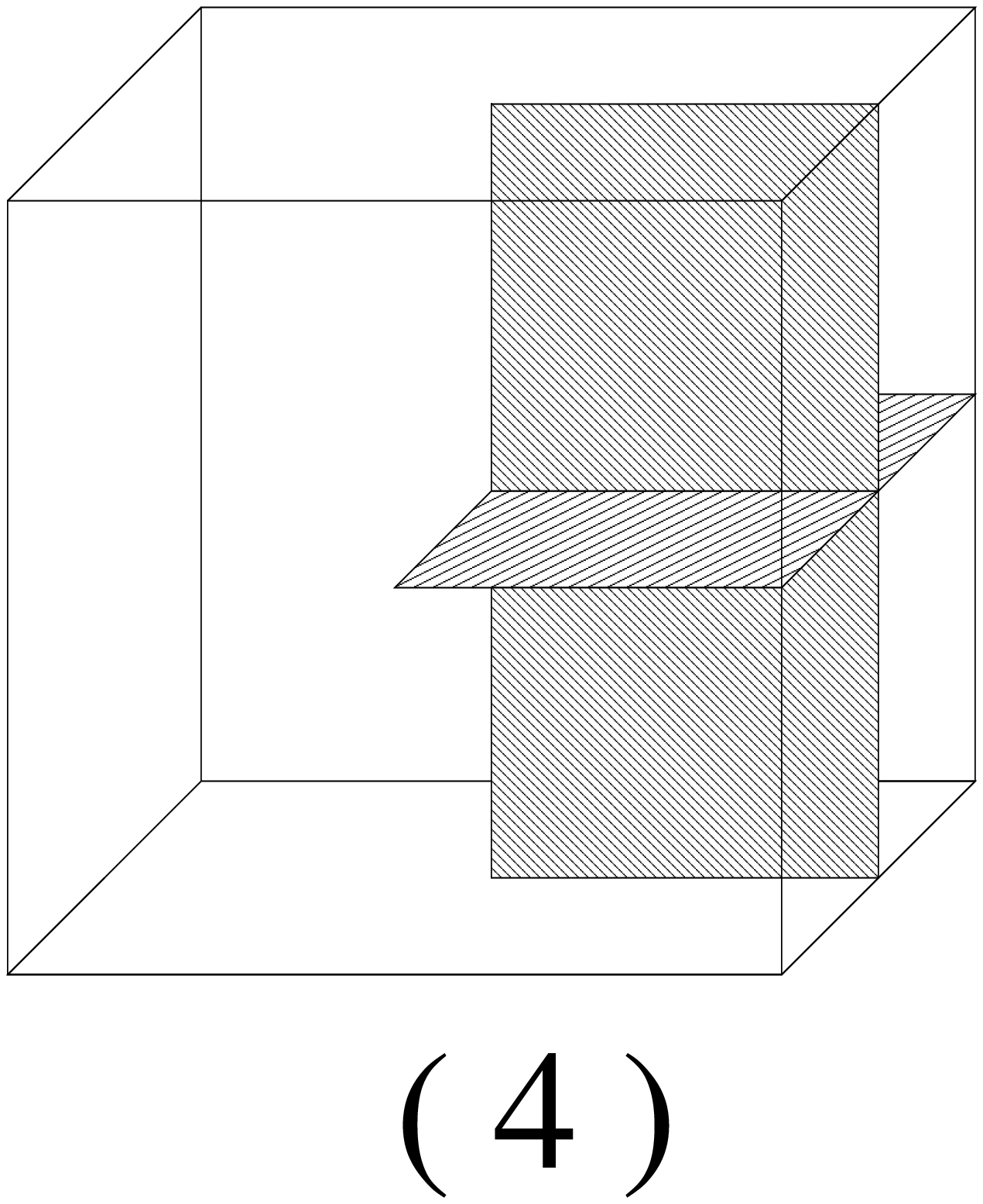}}
\scalebox{.15}{\includegraphics{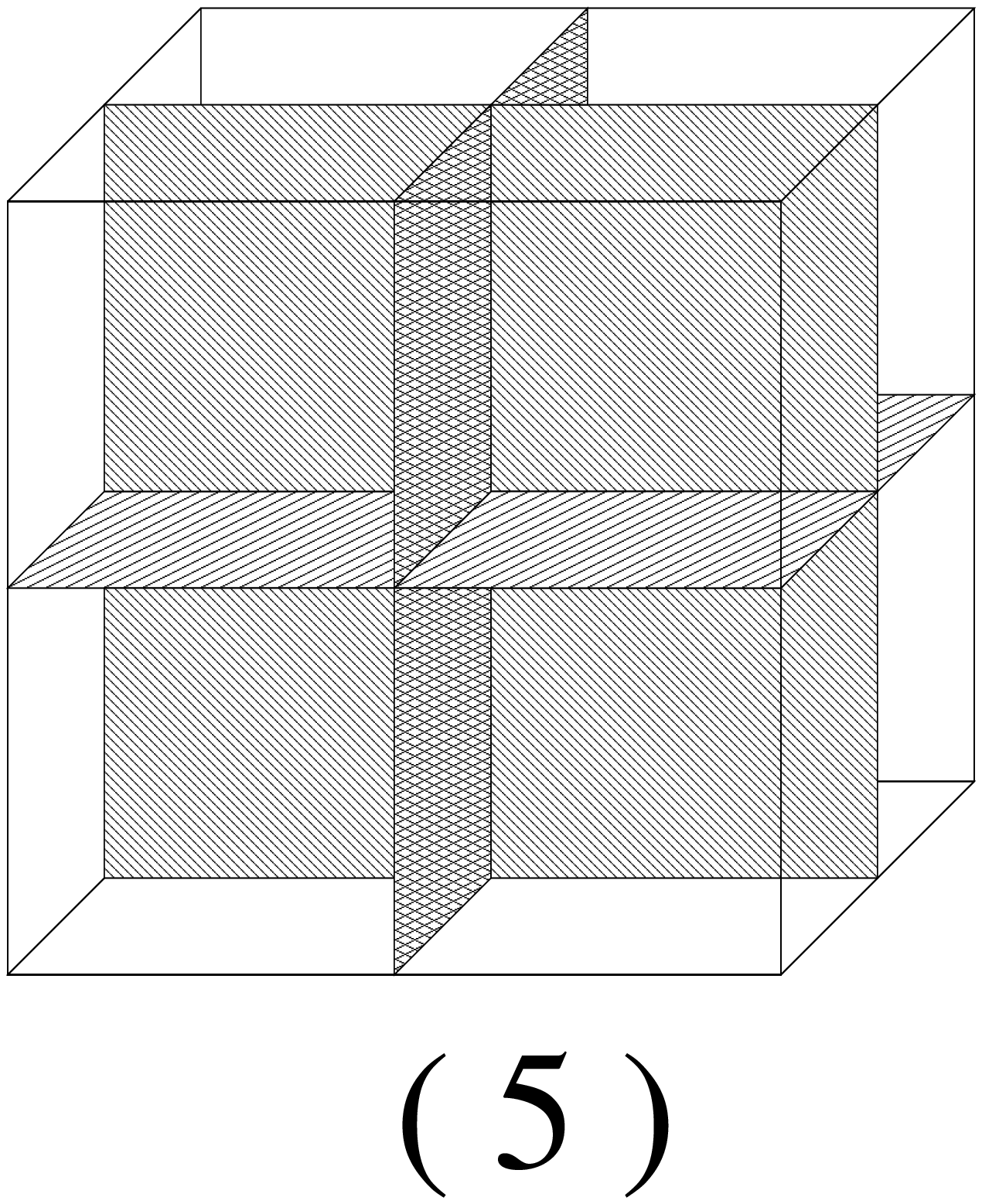}}
\scalebox{.15}{\includegraphics{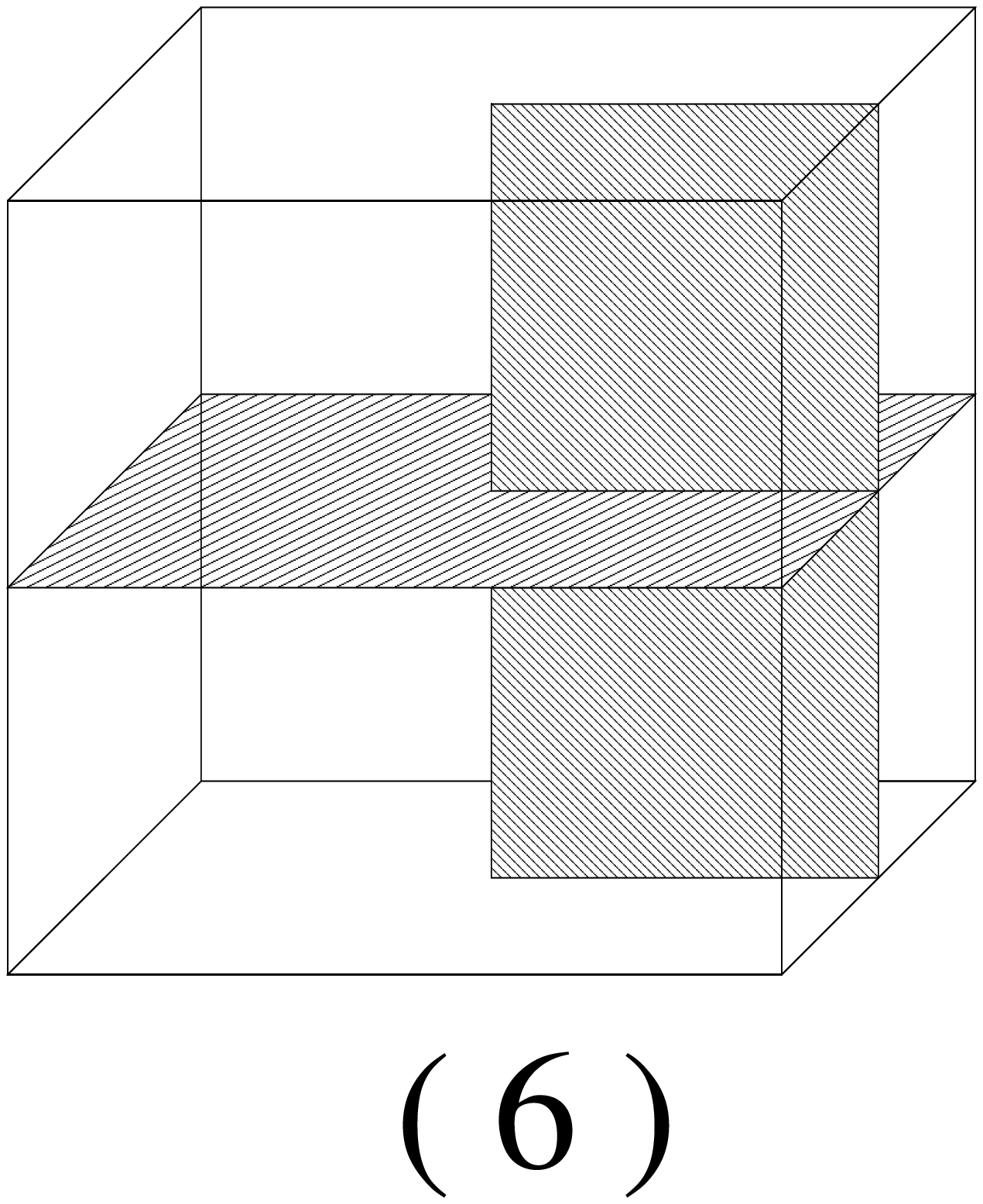}}
\scalebox{.15}{\includegraphics{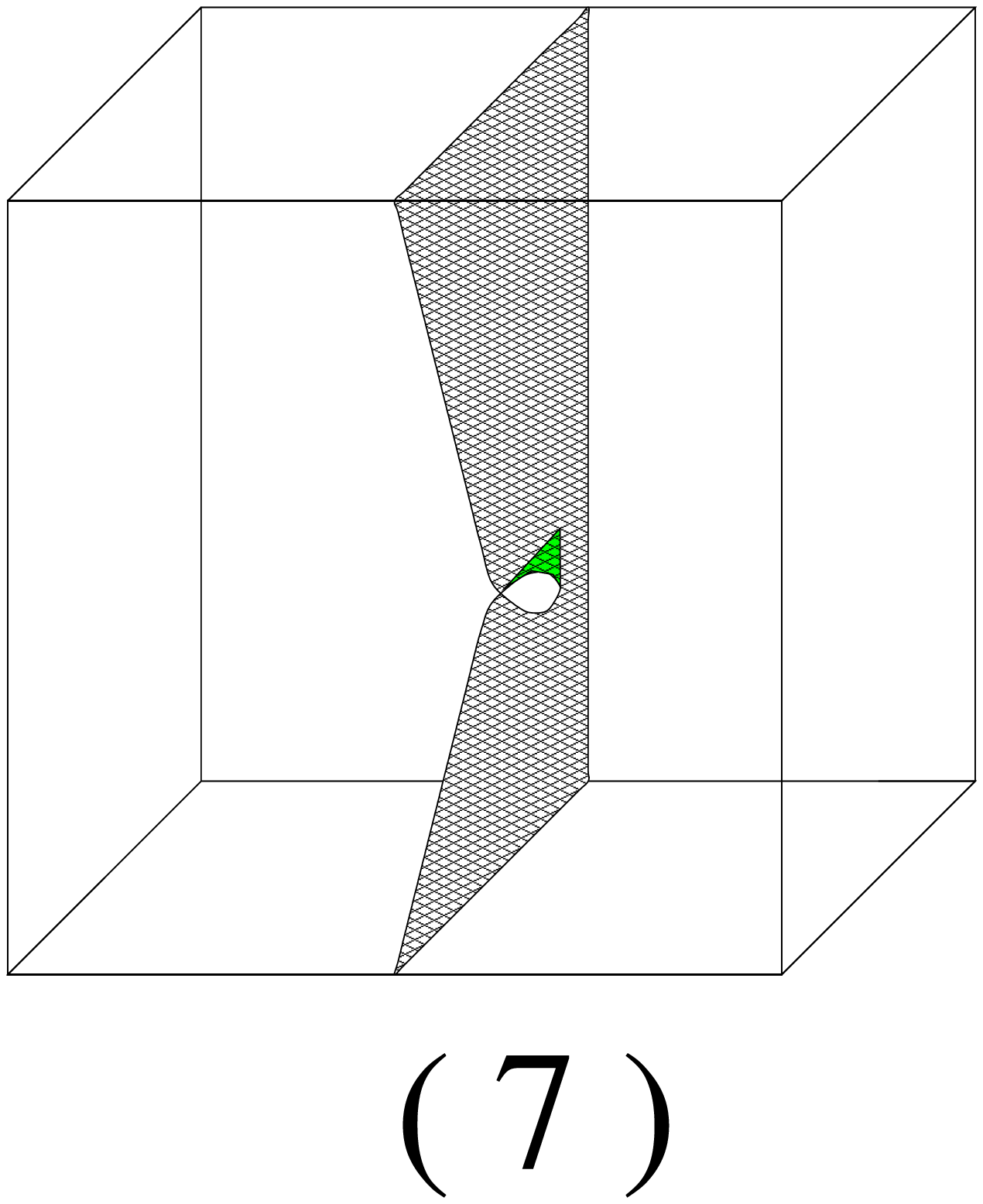}}
\caption{Types of 2-complexes which appear locally in $\tss$} \label{types}
\end{center}
\end{figure}

\begin{thm}{$\pi_2(\tck(\R^3))=\Z$.} \label{pi2}
\end{thm}
\pf \qquad
Let $\gamma: S^2\to\tck$.  Up to homotopy, $\gamma(S^2)\subset\tck^{(2)}$ and is incident to finitely many 2-cells $B_i$.  We can define $h_i: \{x\in S^2|\; \gamma(x)\in B_i\} \to I^2$ and $\Phi: S^1 \times S^2 \to \R^3$ such that $\gamma(x)=g_{\Phi(S^1,x_i)}(h_i(x))$ where $\Phi(S^1,x)$ is a 2-singular knot for finitely many $x_i \in S^2$, and is 1-singular or nonsingular otherwise.  The technical details are similar to those in the proof of Corollary \ref{pi1}, so here we simply regard $\Phi$ as a 2-parameter continuous family of maps $S^1\to\R^3$ which determines a 2-sphere in $CK^{(2)}$.  Thus,
\[ \ss= \text{ closure}\{x \in S^2| \; \Phi(\cdot,x):S^1\to\R^3 \text{ is not an embedding}\} \]
is a graph of valence 1 or 4 on $S^2$.  As above, vertices of valence 1 are limits of maps corresponding to cusps, and vertices of valence 4 are 2-singular knots.  Unlike the case for $S^3$, there is no obstruction to extend $\Phi$ to a map $\tphi: S^1 \times D^3 \to \R^3$.  Using the same almost general position argument as in Theorem \ref{lin}, we can perturb $\tphi$ relative to $\del D^3$ such that
\begin{itemize}
\item[$i)$] $\tss= \text{ closure}\{x \in D^3| \; \tphi(\cdot,x):S^1\to\R^3 \text{ is not an embedding}\}$
is a complex of dimension 2 in $D^3$ with $\tss|_{\del D^3}=\ss$, and
\item[$ii)$] locally, $\tss$ is one of the seven types of 2-complexes shown in Figure \ref{types}.
\end{itemize}

The following definitions follow the proof of Corollary \ref{pi1}.  We define a
\emph{subdivision of $D^3$} to be pairs $\{(R_\alpha,h_\alpha)\}$ such that $D^3=\cup
R_\alpha$ and $h_\alpha: R_\alpha\to I^3$. We will say that a nonempty 2-complex $\tss$ in
$D^3$ \emph{can be realized in $\tck$} if there exists a map $\Psi:D^3\to\tck$ and a
subdivision of $D^3$, such that for $x\in R_\alpha$, $\Psi(x)=g_{\tphi(S^1,x_\alpha)}
(h_\alpha(x))$ for some $x_\alpha\in R_\alpha\cap\tss$.

The 2-complexes of type $(1)$ and $(2)$ can be realized in $CK$.  By extending the proof of
Corollary \ref{pi1}, type $(3)$ can be realized in $CK\cup_fA_1$, and hence in $\tck$.  We now consider types $(4)-(7)$.

\begin{figure}
\begin{center}
\scalebox{.45}{\includegraphics{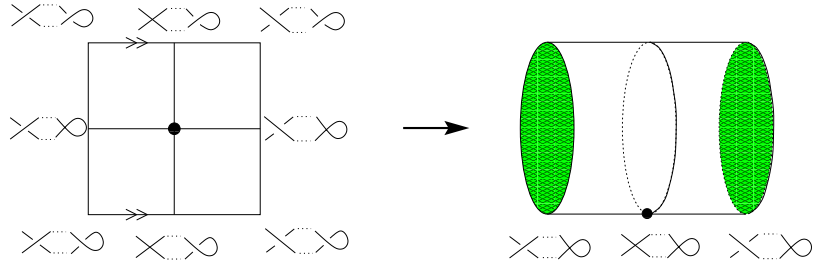}} \caption{2-sphere which appears in type
$(6)$} \label{type6}
\end{center}
\end{figure}
\begin{figure}
\begin{center}
\scalebox{.45}{\includegraphics{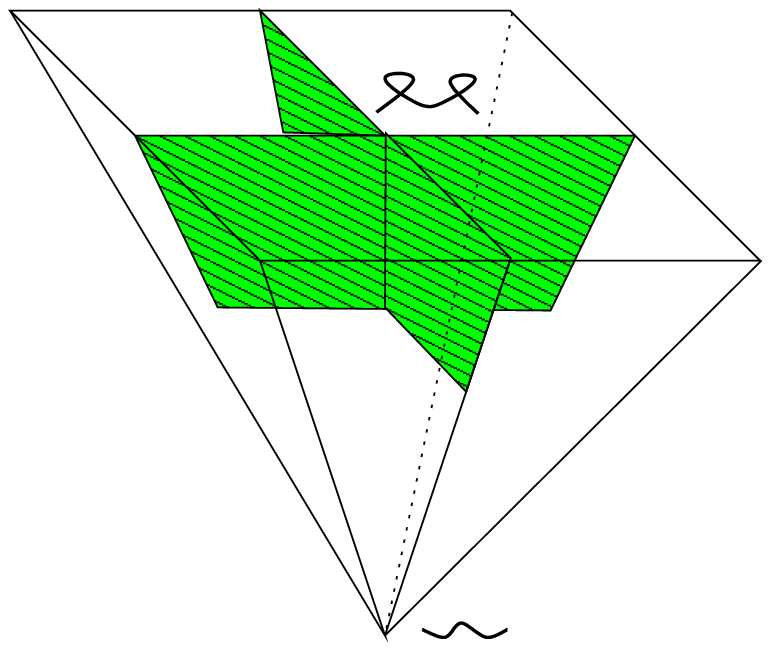}}
\scalebox{.45}{\includegraphics{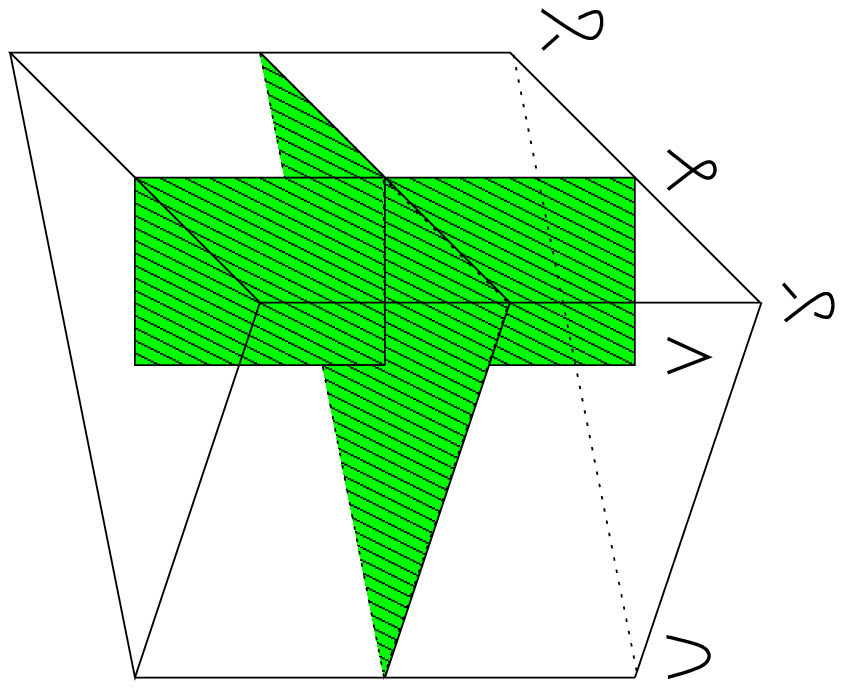}} \caption{Another view of types $(4)$ and $(6)$} \label{type4and6}
\end{center}
\end{figure}

Types $(4)$ and $(6)$ correspond to degenerate 2-singular knots.  For any $K\in D_2, \;
g_K(I^2)$ contains a 2-sphere, as shown in Figure \ref{type6}.  In our construction of $\tck$,
the map $f_{\alpha_2}$ is given precisely by attaching a 3-cell to every such 2-sphere.  By
attaching these cells, types $(4)$ and $(6)$ can be realized in $\tck$.  These operations
resemble some kind of ``compactification'' of $CK$ with respect to its 2-skeleton, as
suggested in Figure \ref{type4and6} where we redraw types $(4)$ and $(6)$.

\begin{figure}
\begin{center}
\scalebox{.35}{\includegraphics{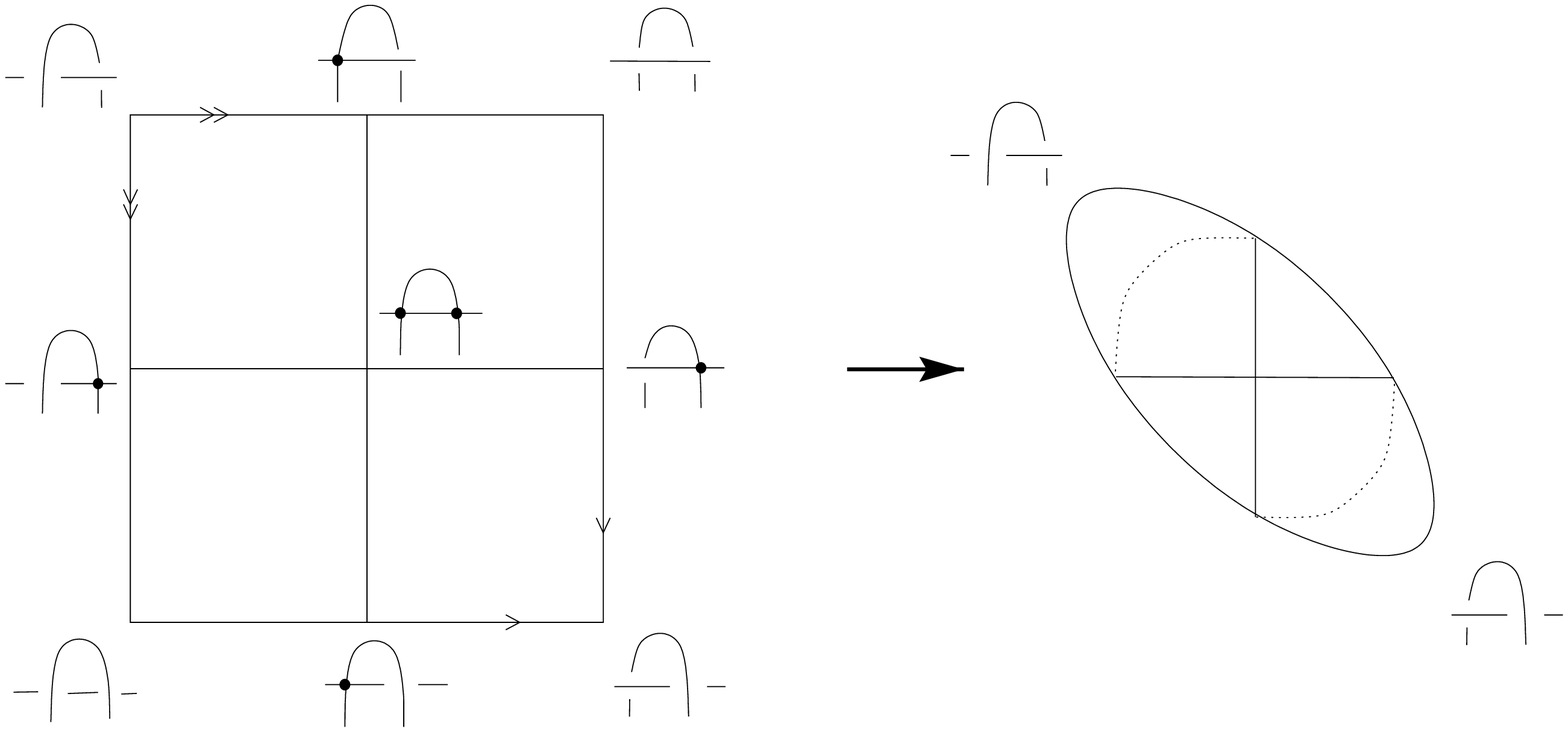}}
\scalebox{.45}{\includegraphics{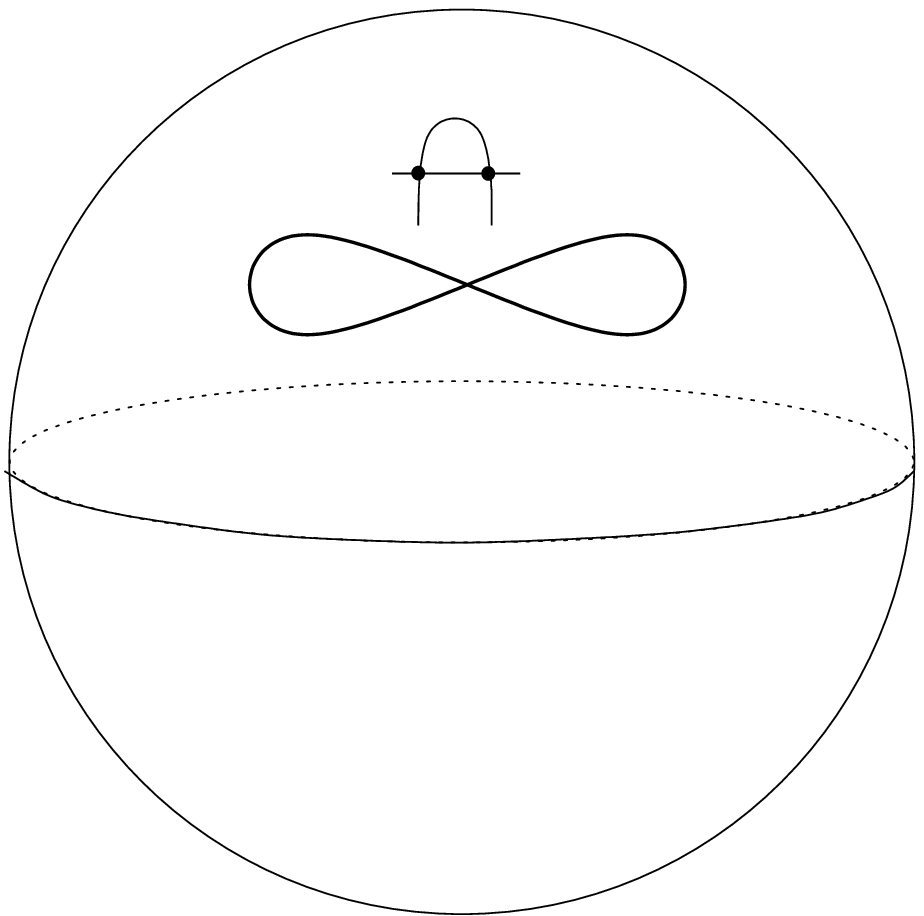}} \caption{Basic model for the boundary
of type $(7)$ and related graph $\ss\subset S^2$} \label{type7}
\end{center}
\end{figure}

The basic model for the boundary of type $(7)$ is the 2-cell shown in Figure \ref{type7}.
The identifications on this 2-cell follow from the topological move shown in Figure
\ref{akirafig}, and therefore its image in $CK$ is a 2-sphere.  Also shown in Figure
\ref{type7} is the corresponding graph $\ss$ on the 2-sphere.

In the basic model for the boundary of type $(7)$, we have shown only fragments of a knot diagram.  We claim that changing any other nonsingular crossings can be realized by a homotopy of 2-spheres in $CK$.
Let $K$ be any 2-singular knot given by $\Phi(S^1,x)$, such that $x$ is the 4-valent vertex in the boundary of type $(7)$.  If $K'$ is obtained from $K$ by changing any nonsingular crossing, then we can find a 3-cell $g_{K_{123}}(I^3)$ in $CK$ with the front face of $I^3$ mapped by $g_K(I^2)$ and back face mapped by $g_{K'}(I^2)$.  We may continue changing crossings, gluing each 3-cell front to back, and keeping the identification of the sides of $I^3$.  The 3-cell obtained in this way is a homotopy in $CK$ from the top 2-sphere $g_K(I^2)$ to the bottom 2-sphere, which can be chosen to be any representative of our basic model.

\begin{figure}
\begin{center}
\scalebox{.3}{\includegraphics{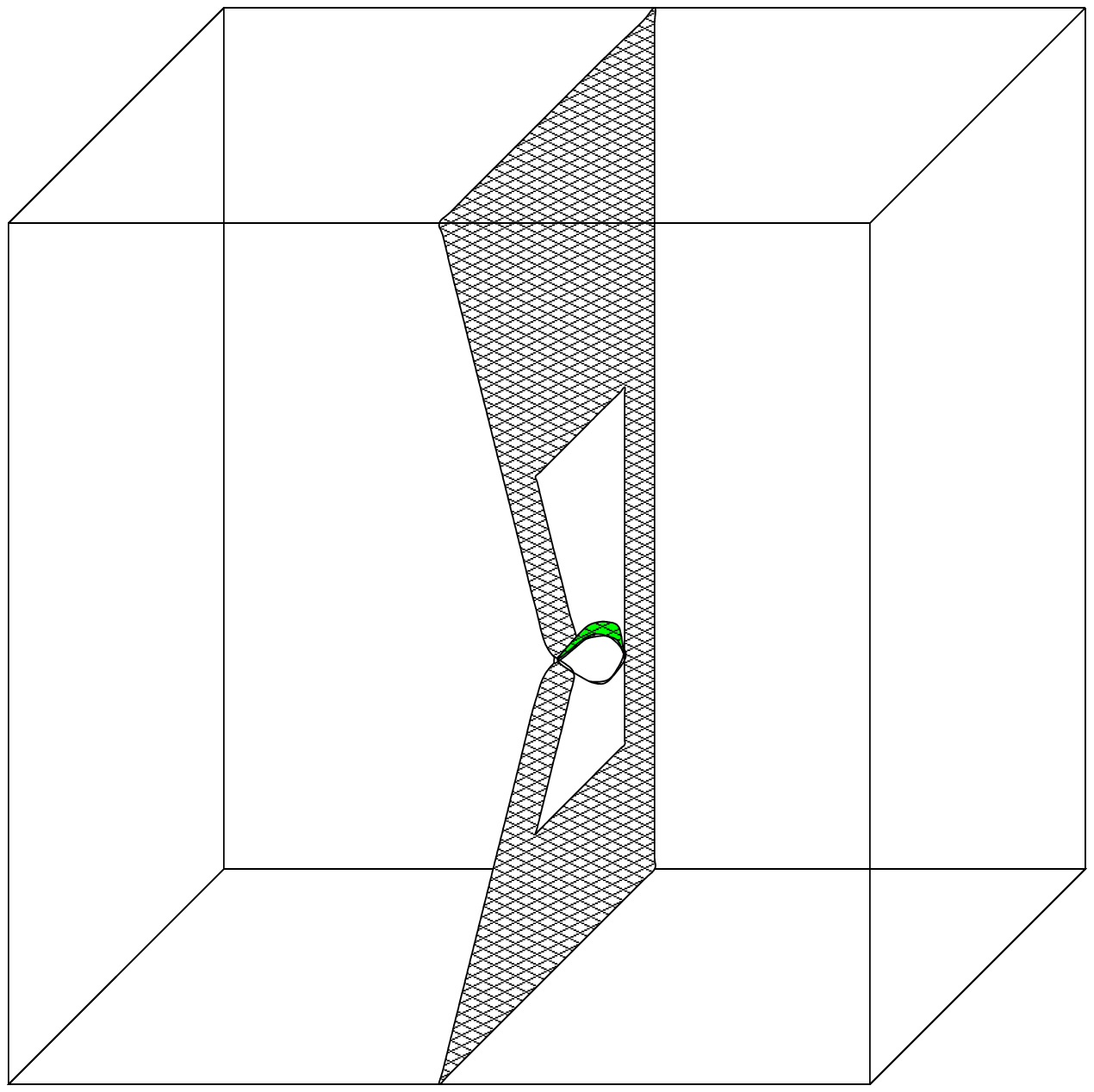}}
\scalebox{.3}{\includegraphics{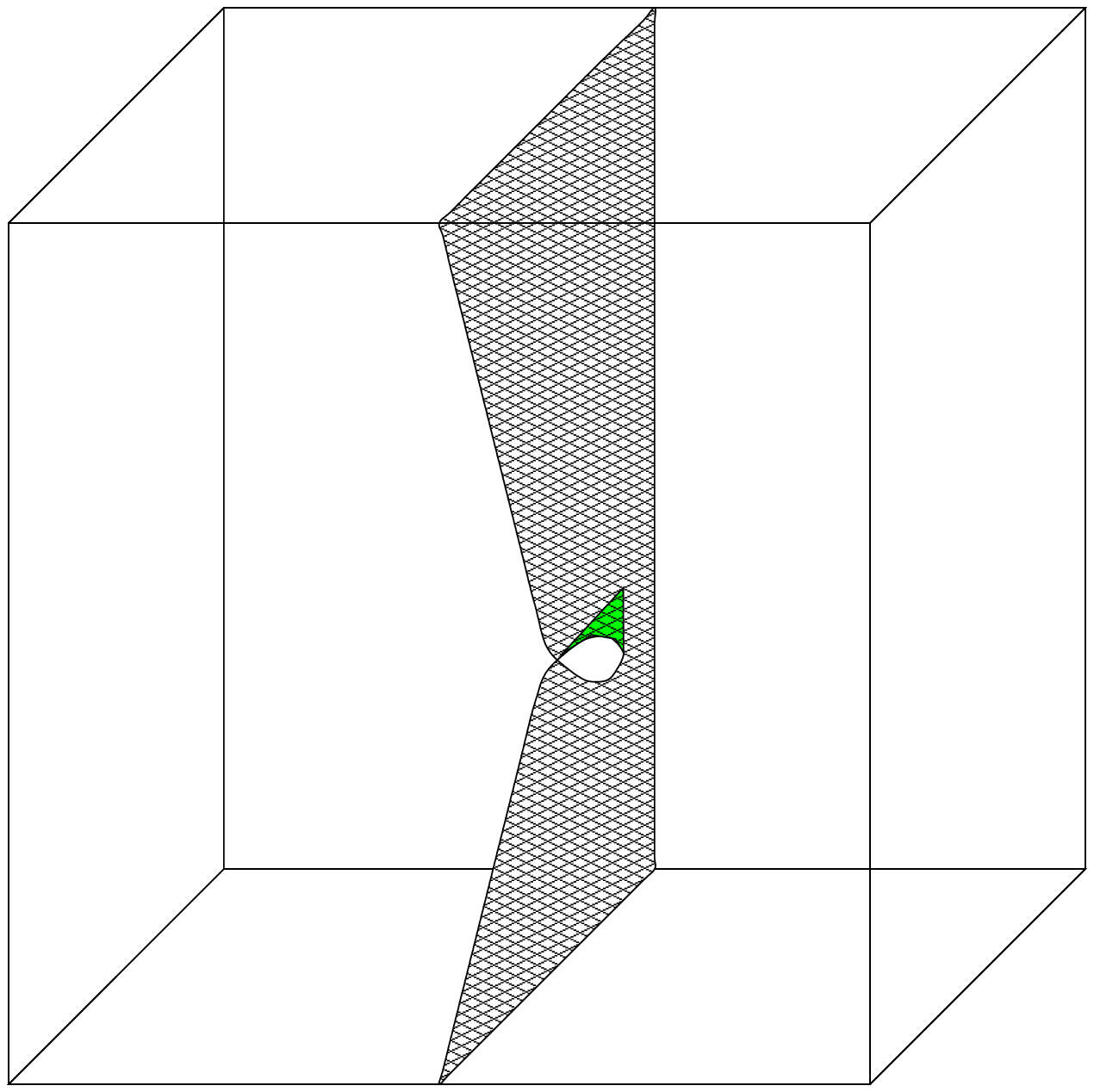}}
\\
\scalebox{.16}{\includegraphics{trivgen.eps}}\hspace*{1in}
\scalebox{.12}{\includegraphics{trefgen.eps}}
\caption{Two ways to extend the basic model from Figure \ref{type7} to a 2-complex, and associated 2-singular knots} \label{type7b}
\end{center}
\end{figure}

Since, up to homotopy, we can change nonsingular crossings, the basic model can be realized as
a 2-sphere in $CK$ in two different ways:
$\mbox{\epsfig{file=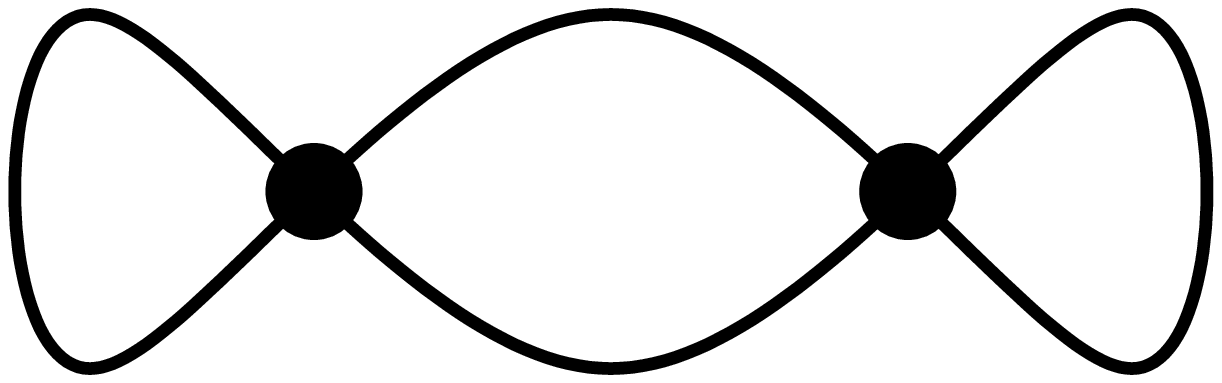,height=11pt}}$ or
$\mbox{\epsfig{file=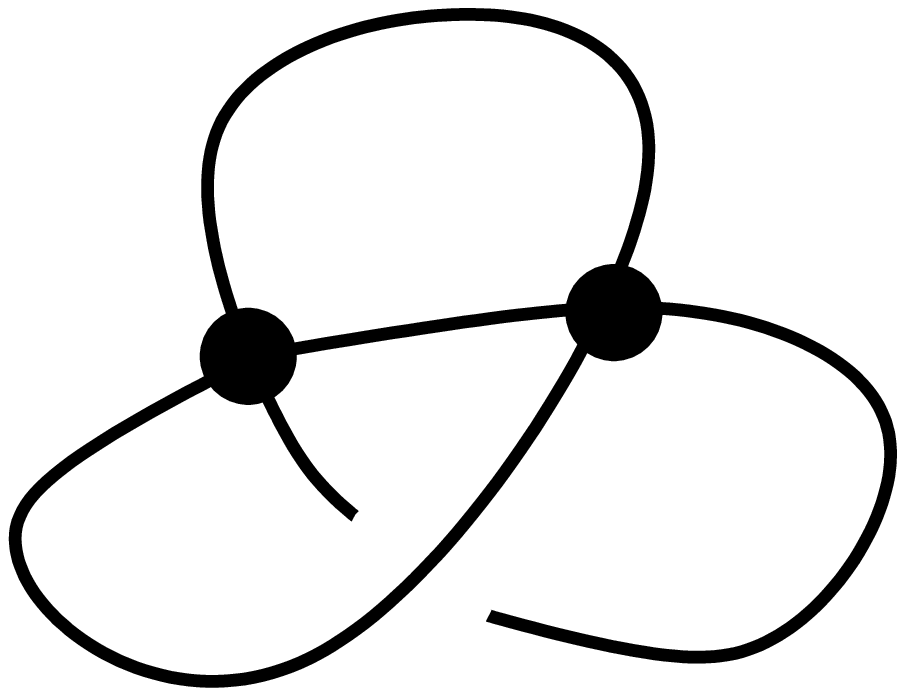,height=15pt}}$.  Up to crossing changes, these are the
two distinct classes of 2-singular knots.  These two cases correspond to the two possible
extensions of the boundary of type $(7)$ to locally different 2-complexes, as shown in
Figure \ref{type7b}.  The associated 2-singular knots (up to crossing changes) are shown below
each 2-complex.
The first 2-complex in Figure \ref{type7b} can be realized in $\tck$ because it is essentially
of type $(4)$.  Thus, the 2-sphere corresponding to
$\mbox{\epsfig{file=trivgen.eps,height=11pt}}$ is trivial in $\pi_2(\tck)$.  

However, the 2-sphere corresponding to $\mbox{\epsfig{file=trefgen.eps,height=15pt}}$ is nontrivial in $\pi_2(\tck)$, and any choice of ordering or base point is equivalent up to homotopy.  The latter claim is that $\mbox{\epsfig{file=trefgen.eps,height=15pt}}$ represents a unique element up to crossing changes in $X_2^b$ and in $X_2^0$, which can be shown by an explicit isotopy, but is obvious by considering its chord diagram.  We can prove nontriviality by using results from Section \ref{sec5}.  For $\Phi:S^1\times I^2\to\R^3$ as shown in Figure \ref{type7}, let $\Psi:S^2\to\tck$ be the
corresponding element in $\pi_2(\tck)$.  By the Hurewicz theorem, Corollary \ref{pi1} implies
that $[\Psi]$ is nontrivial in homotopy whenever it is nontrivial in homology.  By Corollary
\ref{V2tilde}, the fact that $v_2([\Psi])\neq 0$ proves our claim.

The only remaining case is type $(5)$.  The triple point in this 2-complex corresponds to an immersion with 3 parameters to obtain 2-singular knots.  The partitions $3=1+1+1,\; 3=2+1,\; 3=1+2$ give the possible ways of grouping the parameters.  Namely, the first case is that of 3 separate double points, such that each double point can be independently perturbed to obtain a 2-singular knot.  Of course, this case can be realized in $CK$.

\begin{figure}
\begin{center}
\scalebox{.45}{\includegraphics{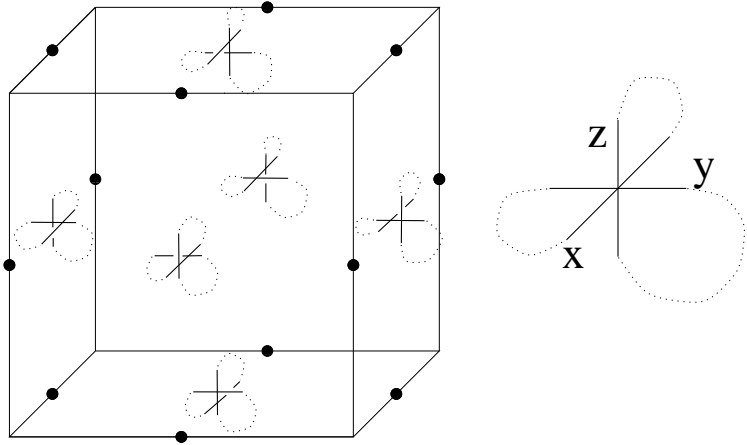}} \caption{$3=2+1$ model for type $(5)$} \label{triple}
\end{center}
\end{figure}
The other two cases, equivalent by symmetry, arise from a triple point on a singular knot.  As
shown in Figure \ref{triple}, the model is to perturb the $z$-strand by 2 parameters, and to
perturb the $x$-strand by 1 parameter to obtain 2-singular knots.  (A similar figure appears
as Figure 36 in \c{V}.)  The boundary of the cube consists of six faces, corresponding to any
realization of the six 2-singular diagram fragments shown, and can be mapped to a 2-sphere in
$CK$.  By a similar argument to that above, changing any nonsingular crossings can be realized
by a homotopy along 3-cells in $CK$. Therfore, the boundary of the cube is homotopic to the
boundary of the cube formed by $3\,\mbox{\epsfig{file=trefgen.eps,height=15pt}}$ and
$3\,\mbox{\epsfig{file=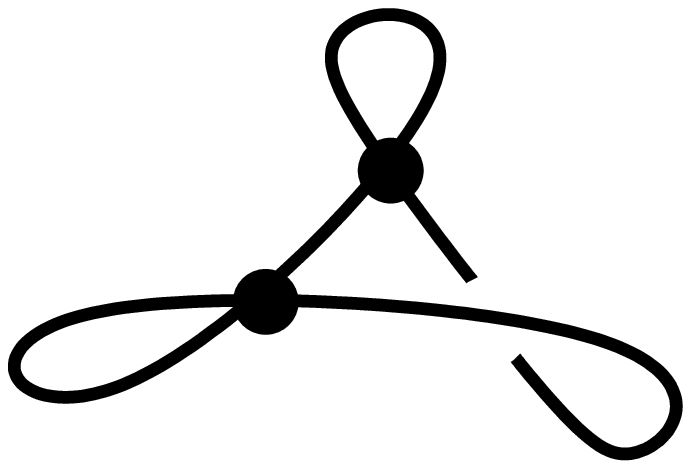,height=15pt}}$.  Since the latter is trivial in
homotopy, the boundary of the cube is a multiple of
$\mbox{\epsfig{file=trefgen.eps,height=15pt}}$.

Therefore, $\pi_2(\tck)\cong \Z$, generated by $\mbox{\epsfig{file=trefgen.eps,height=15pt}}$.
\eop

\section{Comparison with the Vassiliev construction} \label{mark}

In this section, we describe a subcomplex of $CK_b(\S)$ which can be associated to the space $\s$ defined by Vassiliev \c{V}.

Instead of knots in $\S$, Vassiliev considered the equivalent space $\mathscr{M}$ of embeddings $\R^1 \to \R^3$ that asymptotically approach a fixed line (often called \emph{long knots}).  These can be identified with elements of $X_0^b(\S)$, with the base point placed at infinity.  Let $\Si$ denote the discriminant of this space; i.e., the set of maps with singularities or self-intersections.  Let $\Gamma^d \subset \mathscr{M}$ be the finite-dimensional subset generated by certain monic polynomials of degree $d+1$, suitably approximated to be transverse to $\Si$.  The homology of $\s$ coincides with that of $\Gamma^d \cap \Si$.

The nonsingular part of $\Si$ is the set of immersions with exactly one transverse double point.  A nonsingular self-intersection point of $\Si$ is an immersion with exactly two distinct double points.  Similarly, nonsingular \n-fold self-intersection points of $\Si$ correspond to immersions with $n$ distinct double points.  The essential idea of Vassiliev invariants is that once their value is defined on the unknot, their value on any other knot is equal to the algebraic sum of indices assigned to these points of self-intersection, as if passing from one component to another and recording the ``jump'' of the invariant at the walls.

We can encode double points of the knot by specifying pairs of points on the line $\R^1$, and triple points of the knot by three specified points on the line.  To construct $\s$, Vassiliev showed that it was sufficient to consider \emph{noncomplicated} cases, meaning only double points, or double points and exactly one triple point, or double points and exactly one other singular point.  In fact, it suffices to consider only the first two cases (see \c{BL}).  Thus, we can simplify matters and let an \emph{$i$-configuration} denote the collection of points on the line (equivalent up to orientation-preserving diffeomorphism of $\R^1$) corresponding to $i$ double points, or $i-2$ double points and a triple point.  We say a map $\Theta: \R^1 \to \R^3$ respects an $i$-configuration if it has only these specified singularities.  Let $\mm(\Gamma^d,J)$ be the set of maps in $\Gamma^d$ which respect $J$.

To any $i$-configuration $J$ and any map $\Theta \in \mm(\Gamma^d,J)$, Vassiliev defined a
corresponding simplicial complex in $\R^N$.  The first step is to define \emph{generating
collections} for $J$ by perturbing the triple point (if there is one) to two double points, as
in Figure \ref{triple} (Figure 36 \c{V}).  To any generating collection with $\ell$ double
points, he associated a specific $(\ell-1)$-dimensional simplex in $\Theta \times \R^N$ (for
details, see section $V.2.3$ \c{V}).  Let $\s \subset \Gamma^d \times \R^N$ be the union of
simplices over all $i$-configurations, for $i \leq 3d$.

As pointed out in \c{BL}, if instead we let $S_J$ be the simplex associated to the unique maximal generating collection of $J$, then all the other simplices associated to $J$ are sub-simplices of $S_J$.  We therefore obtain:
\[ \s=\bigcup \mm(\Gamma^d,J) \times S_J \subset (\Gamma^d \cap \Sigma) \times \R^N \]
where the union is taken over all $i$-configurations $J$, for $i \leq 3d$.  If we consider the projection on the second factor $p:\s \to \R^N$, the image $p(\s)$ is a simplicial complex.  For any point $x$ in the interior of a simplex of $p(\s)$, the preimage $p^{-1}(x)$ is the affine space $\mm(\Gamma^d,J)$ \c{BL}.

In addition, there is a filtration $\s_1 \subset \cdots \subset \s_{3d} = \s$, where $\s_i \setminus \s_{i-1}$ is the union of simplices corresponding to $i$-configurations.  We obtain the \emph{canonical decomposition} of $\s_i \setminus \s_{i-1}$ by decomposing the second factor into open cells corresponding to various generating collections of $J$.  The canonical decomposition together with an added point define a CW-complex structure on the one-point compactification of $\s_i \setminus \s_{i-1}$ (Proposition 3.2.2 \c{V}).


\begin{figure}
\begin{center}
\scalebox{.45}{\includegraphics{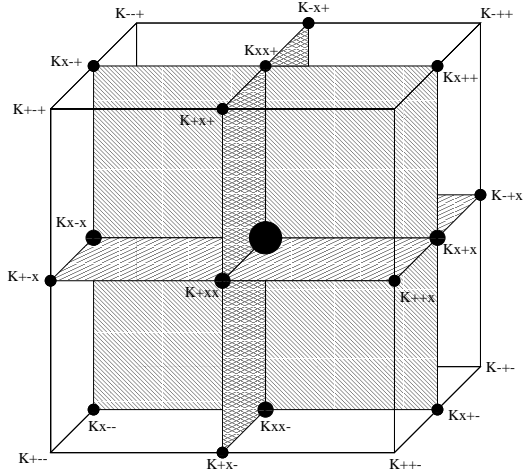}}
\caption{Markings on the 3-cube labeled by $K_{123}$}
\label{markedKxxx}
\end{center}
\end{figure}

We now describe a subcomplex of $CK_b(\S)$ which can be associated to the discriminant $\Si$. Recall $I^n=\{\vec{x} \in \R^n : -1 \leq x_i \leq 1 \}.$  Let $\Si(I^n)$ denote the hyperplanes $\{x_i=0\}, 1\leq i \leq n$, which are called \emph{markings}.  Let $\Si_k(I^n)$ denote the $k$-skeleton of $\Si(I^n)$, subdivided into $\binom{n}{k} \cdot 2^k \; k$-planes, each incident to $\vec{x}=0$.  $\Si_0(I^n)=\{0\}$.  (See Figure \ref{markedKxxx}.)

\begin{defn}
For any $K \in X_n^b(\S)$, let $\Si(K)=g_{K}|_{\Si(I^n)}$ and let
$\Si_k(K)=g_{K}|_{\Si_k(I^n)}$.  Let $\Si^n(CK_b)=\cup_{X_n^b}\Si(K)$.
\end{defn}

For any $K \in X_n^b(\S)$, with the base point placed at infinity, let $\pi(K)$ be the projection to the corresponding \n configuration.  Let $\mm_{K}$ be the rigid-vertex isotopy class of maps with distinct double points and no other singularities, considered as a proper subset of $\bigcup_{d \geq 0}\mm(\Gamma^d, \pi(K))$.  Fix $d>>n$.  For each $\Theta \in \mm_{K} \cap \Gamma^d$, there is an associated simplex $s_{\pi(K)}$ in $\R^N$, as described by Vassiliev.  Then, $\Theta \times s_{\pi(K)}$ is contained in a cell of the canonical decomposition of $\s_n \setminus \s_{n-1}$. 

We now define the map $\Psi^n: \s_n\to\Si^n(CK_b)$.  For any $\Kn\in X_n^b(\S)$, to define the
corresponding attaching map for $I^n$, we defined a labeling on $\del I^n$, as in Figure
\ref{Kxxx}.  Let $K_{n-k}$ be any $(n-k)$-singular knot which is a resolution of $\Kn$.  For
any $k,\; 0 \leq k \leq n-1,$ there are $\binom{n}{k} \cdot 2^k$ such resolutions and
corresponding labels on $\del I^n$.  (See Figure \ref{markedKxxx}.)  Let $\Psi^n_0(\mm_{\Kn}
\times s_{\pi(\Kn)}) = \Si_0(\Kn)$.  For any $k,\; 0 < k \leq n-1$, let $\rho:
s_{\pi(K_{n-k})}\to\Si_k(I^n)$ be the simplicial projection to the corresponding $k$-plane,
according to the labeling for $\Kn$.  Let $\Psi^n_k:\mm_{K_{n-k}} \times s_{\pi(K_{n-k})} \to
\Si(\Kn)$ be defined by projecting to the second factor and applying $g_{\Kn}\circ\rho$.  This
gives an identification between $\Si_k(\Kn)$ and $(n-k)$-singular knots which are resolutions
of $\Kn$:  for any point $x$ in the interior of $\Si_k(\Kn)$ corresponding to the resolution
$K_{n-k}$, the preimage $(\Psi^n_k)^{-1}(x)$ is the space $\mm_{K_{n-k}}$.

As mentioned above, for $i\geq 2$, an $i$-configuration may also contain $i-2$ double points
and one triple point, where the triple point is resolved as in Figure \ref{triple} (Figure 36
\c{V}). For $n\geq 0$, let $K_{(n,1)}$ denote an \n singular knot with an additional triple
point apart from the double points.
\begin{conj}
Let $CK^*$ be the space obtained from $CK_b(\S)$ by attaching an $(n+1)$-cell for every knot $K_{(n-2,1)},\; n\geq 2$, with boundary given by resolving the double points,
and resolving the triple point as in Figure \ref{triple}.  Then, $\bigcup_n \Si^n(CK^*)$ is
homotopy equivalent to the Vassiliev simplicial complex.
\end{conj}


\section{Vassiliev invariants and cohomology of $CK$} \label{sec5}

A knot invariant is a map $\phi: X_{0} \to \Q$.  Any knot invariant can be extended to singular knots by inductively using the following skein relation:

\begin{equation} \label{skein}
 \phi(\Kx) = \phi(\Kp) - \phi(\Km), \quad \phi(K_{\times \ldots \times \ldots \times}) = \phi(K_{\times \ldots + \ldots \times}) - \phi(K_{\times \ldots - \ldots \times})
\end{equation}

\begin{defn} \label{v}
A knot invariant is said to be of \emph{finite type} or a \emph{Vassiliev invariant} if there exists $n \in \mathbf{N}$ such that $\phi(K_{\underbrace{\xdots}_{n+1}})=0$.  An invariant of singular knots is of \emph{order n} if $n$ is the smallest such integer.  We denote by $\V_n(\M)$ the vector space of invariants of knots in $\M$ of order $\leq n$.
\end{defn}

We can extend any singular knot invariant $\phi: X_n \to \Q$ to a cochain $\phi \in C^n(CK)$ by forgetting the ordering and/or base point: $\phi(\Kn)= \phi(\Kdots)$.  In particular, the vector space over $\Q$ freely generated by restricting $\V_N|_{C_n}$ can be made a subspace of the cochain complex $\{C^n(CK,D), \delta^n\}$ in this way.  In addition, for any Vassiliev invariants $v,w \in \V_N$ if $v|_{C_n}=w|_{C_n}$ then on the space of knots $v=w$ up to invariants of order $n-1$.  Therefore, $\V_N/\V_{n-1} \hookrightarrow C^n(CK,D)$ is an embedding.  We adopt the following standard notation:
\[ \begin{array}{l}
Z^n(CK,D)= ker (\delta^n), \; n \text{-cocycles}, \\ \\
B^n(CK,D)= image (\delta^{n-1}), \; n \text{-coboundaries}, \\ \\
H^n(CK,D)= Z^n(CK,D)/B^n(CK,D), \; \text{cohomology with coefficients in } \Q. \\
\end{array} \]

\begin{prop}\label{zn}  For any $\phi \in \V_n, \; \delta^n \phi=0.$ \; i.e.,
$\V_n/\V_{n-1} \hookrightarrow Z^n(CK)$.
\end{prop}
\pf\qquad For any $K_{1 \ldots n+1}\in X_{n+1}^b\text{ or } X_{n+1}^0, \; \phi(\di K_{1 \ldots n+1})=0 \quad \forall i$.
\[ \delta \phi(K_{1 \ldots n+1}) = \phi(\del K_{1 \ldots n+1}) = \phi \left( \sum\limits_{i=1}^{n+1} (-1)^{i+1} \di K_{1 \ldots n+1} \right) =0 \qquad \square \]

From the exact sequence of the pair $(CK,D)$, we obtain the following commutative diagram:
\begin{equation} \label{diag}
\xymatrix{ {H^{n-1}(D;\Q)} \ar[r] & {H^n(CK,D;\Q)} \ar[r]^\alpha & {H^n(CK;\Q)} \ar[r]^\beta & {H^n(D;\Q)} \\
{} & {} & {\V_n|_{C_n}} \ar[ul]^{\exists\gamma} \ar[u] \ar[ur]_0 & {} }
\end{equation}
Since $\V_n|_{C_n} \subset \ker \beta = \text{Im} \, \alpha$, there is an induced map $\gamma: \V_* \rightarrow H^*(CK,D)$.  Understanding this map is our main problem in this section.

\subsection{Low-dimensional cohomology}

\begin{thm}{For $\phi \in \V_2(\R^3), \; 0 \neq [\phi] \in H^2(CK(\R^3)).$}\label{V2}
\end{thm}

\pf\qquad
Let $T$ be any nondegenerate singular knot in $X_2(\R^3). \quad \V_2(\R^3)$ is generated
by $\phi$ such that $\phi(T)=1$ and $\phi(K_{\times\times\times})=0$.  If $T' \in X_2$ is obtained from $T$ by crossing changes, then $\phi(T')=\phi(T)$, so we may assume
\[ T = \mbox{\epsfig{file=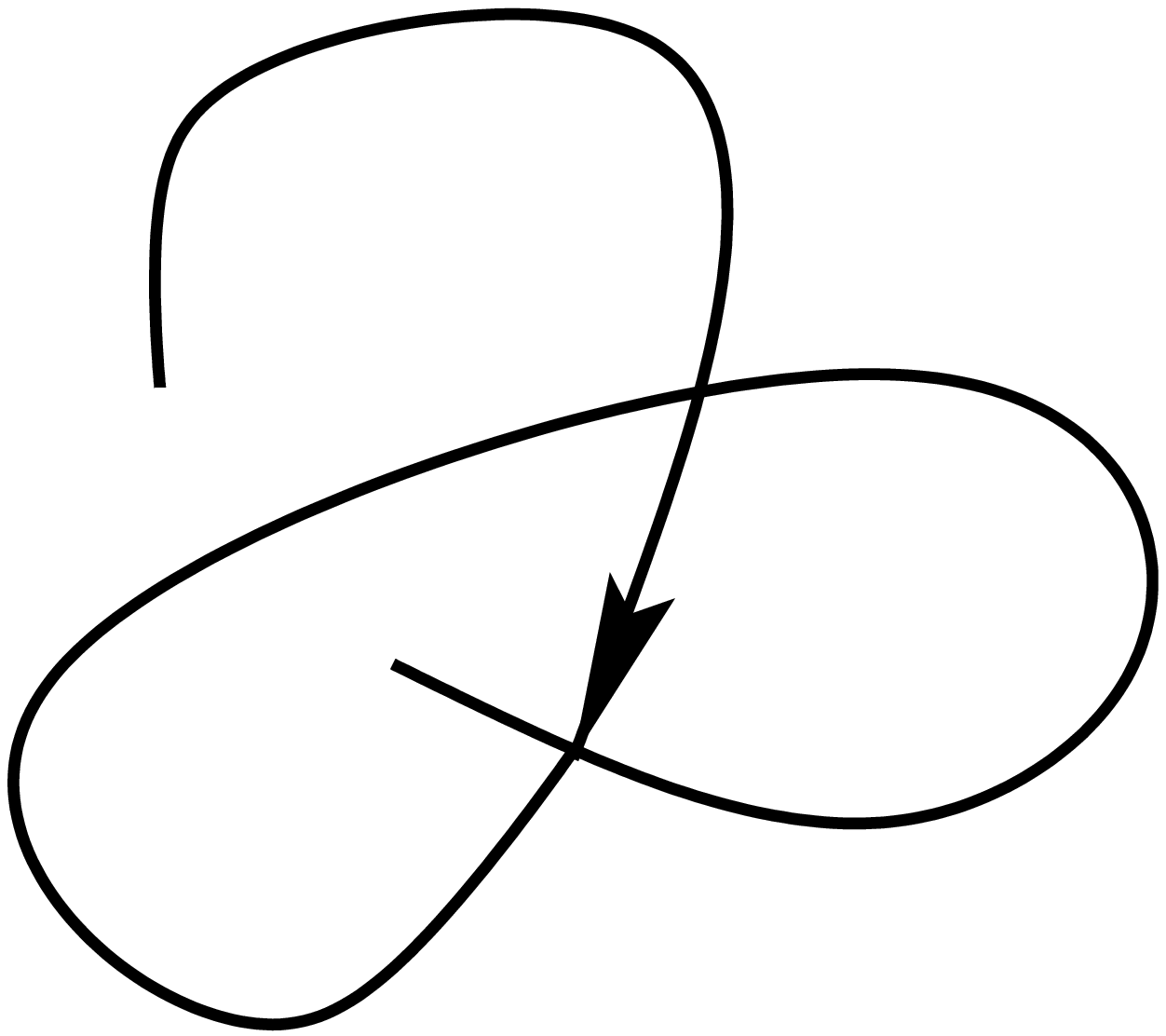,height=15pt}}\; .\]
By Proposition \ref{zn}, $\delta \phi=0$.  Suppose by contradiction $\phi=\delta \psi$ for some $\psi \in C^1(CK)$.  We consider any corresponding element $T\in X_2^0(\R^3)$ or $X_2^b(\R^3)$.
\[ 1=\phi(T)= \delta \psi(T)=\psi(\del(T))= \psi(\d{1}T - \d{2}T) =\psi(0)=0 \quad\Rightarrow\Leftarrow \qquad \square \]

Recall $\tck=CK\cup_f A$, where we attach a 2-cell for every degenerate 1-singular knot, and a 3-cell for every degenerate 2-singular knot.  We can extend any $v \in \V_N$ to $C^*(\tck)$ by $v|_\s=v|_{\del\s}=0, \, \forall \s \in A$.

\begin{cor}{For $\phi \in \V_2(\R^3), \; 0 \neq [\phi] \in H^2(\tck(\R^3)).$}\label{V2tilde}
\end{cor}
\pf\qquad
Since $CK$ is a subcomplex of the CW-complex $\tck$, the inclusion $CK \hookrightarrow \tck$ is a cofibration.  Therefore,
\[ \th^*(\tck,CK)\cong\th^*(\tck/CK)\cong\th^*\left(\bigvee\nolimits_{\alpha\in D_1}S_\alpha^2\bigvee\nolimits_{\beta\in D_2}S_\beta^3 \right) \]
From the exact sequence of the pair $(\tck,CK)$, we obtain the following commutative diagram:
\begin{equation}\label{H2surj}
 \xymatrix{ {\th^2(\bigvee S_\alpha^2 \bigvee S_\beta^3)} \ar[r] & {\th^2(\tck)} \ar[r]^\alpha & {\th^2(CK)} \ar[r]^\delta & {\th^3(\bigvee S_\alpha^2 \bigvee S_\beta^3)} \\ 
{} & {} & {\V_2|_{C_2}} \ar[ul]^{\exists\tilde{\gamma}} \ar[u] \ar[ur]_0 & {} }
\end{equation}
Since $v|_\s=0,\, \forall \s \in A,\, \V_2|_{C_2} \subset \ker \delta = \text{Im} \, \alpha$.  Thus, there is an induced map $\tilde{\gamma}:\V_2|_{C_2}\to \th^2(\tck)$, which is nontrivial by Theorem \ref{V2}.
\eop

\begin{cor} \label{H2bar}
$\V_2(\R^3)|_{C_2}\cong H^2(\tck(\R^3))$.
\end{cor}

\pf\qquad By the Hurewicz theorem, Corollary \ref{pi1} and Theorem \ref{pi2} imply that
$H^2(\tck(\R^3); \Z)=\Z$.   Therefore by Corollary \ref{V2tilde}, $\V_2$ uniquely generates
$H^2(\tck(\R^3))$. \eop

Many arguments using Vassiliev invariants depend on projections of links and graphs, and then
using augmented Reidemeister moves on these projections.  Such arguments are available only in
$\R^3,$ or equivalently in $\S,$ since we can always puncture $\S.$ Theorem \ref{V2} relies on
such arguments, and therefore holds only for $\M=\R^3$ or $\S$.  In \c{L}, these concepts were
extended to more general 3-manifolds by studying Map$(S^1 \times D^2, \M))$ in almost general
position.  In particular, Theorem \ref{integrable} is a generalization of one of the central
results in the theory of Vassiliev invariants.  It was first proved by Stanford \c{S1} for
$\R^3$, and has been extended further by Kalfagianni to closed, oriented, irreducible
3-manifolds which are not ``small'' Euclidean Seifert manifolds (see Theorem 4.1 \c{Kalf}).
Since Theorem \ref{integrable} is important for most results in this section, we assume
$\pi_1(\M)=\pi_2(\M)=1$ in this section.

\begin{defn}
We say $\phi$ is \emph{differentiable} if $\phi(\di \Kn)=\phi(\d{j} \Kn) \quad
\forall i,j$.
\end{defn}
In fact, if we forget the ordering, $\phi: X_n \to \Q$ extends to $X_{n+1}$ satisfying the skein relation (\ref{skein}) if and only if $\phi$ is differentiable (see, e.g., Lemma 3.6 \c{S2}).

\begin{defn}
Let $\phi:X_n \to \Q$.  We say $\phi$ is \emph{integrable} if $\exists\psi:X_{n-1} \to \Q$ such that, for any ordering, $\phi(\Kn)=\psi\circ\di(\Kn) \; \forall i$.  In this case, we write $\int\phi=\psi$.
\end{defn}
We will say $\phi\in C^n(CK)$ is integrable if $\phi$ is an integrable invariant of $X_n$, and is therefore invariant under changes of ordering and/or base point.

\begin{thm}[\c{L}] \label{integrable}
Suppose $\pi_1(\M)=\pi_2(\M)=1.$  Then $\phi:X_n(\M) \to \Q$ is integrable if and only if it satisfies: $(i)$ the $1$-term relation, $(ii)$ the $4$-term relation, $(iii) \; \phi$ is differentiable.
\end{thm}

\begin{cor} \label{C1}
If $\pi_1(\M)=\pi_2(\M)=1$, then $H^1(CK(\M),D)=0.$
\end{cor}
\pf\qquad
If $\phi \in Z^1(CK,D)$, we can find $\psi \in C^0(CK)$ such that $\phi=\delta\psi$.
\begin{itemize}
\item[$(i)$] $\phi|_{D_1} = 0 \Ra \phi$ satisfies the $1$-term relation.
\item[$(ii)$] The $4$-term relation is trivially satisfied for $X_1$.
\item[$(iii)$] $\delta\phi=0 \Ra \phi(\d{1} K_{12})=\phi(\d{2} K_{12})$ so $\phi$ is differentiable.
\end{itemize}
For $CK_0$, Theorem \ref{integrable} implies that $\phi$ is integrable, so $\exists \psi$ with $\phi(\Kx)=\psi(\mathrm{d}\Kx)=\delta\psi(\Kx)$.  For $CK_b$, we must also show that $\phi$ is invariant under changes of base point, so the result follows from the following lemma. \eop

\begin{lem}
Suppose $\pi_1(\M)=1.$  If $\phi\in C^1(CK_b(\M))$ is differentiable, then $\phi$ is invariant under changes of base point.
\end{lem}
\pf\qquad
Since $\phi$ is differentiable, we obtain the following identity for changing crossings and base point:
\[ \begin{array}{l}
\phi(\mbox{\epsfig{file=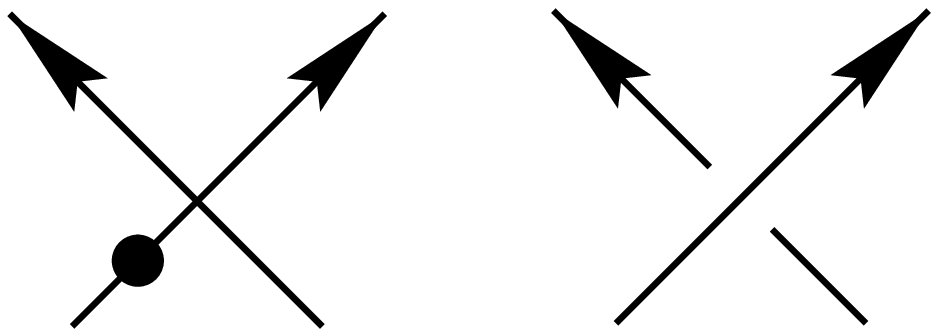,height=15pt}}) - \phi(\mbox{\epsfig{file=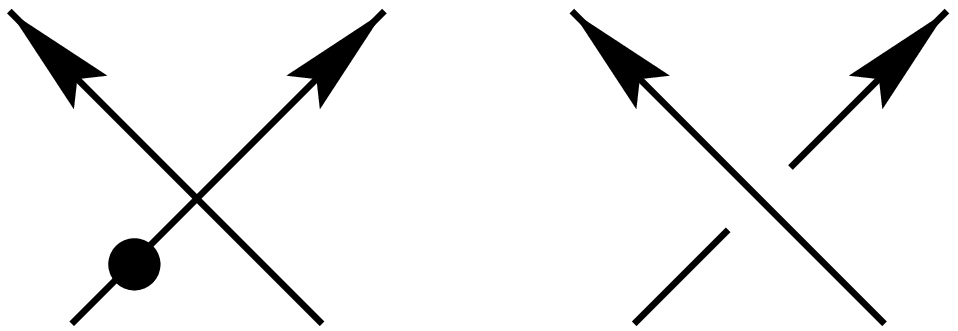,height=15pt}}) =
\phi(\mbox{\epsfig{file=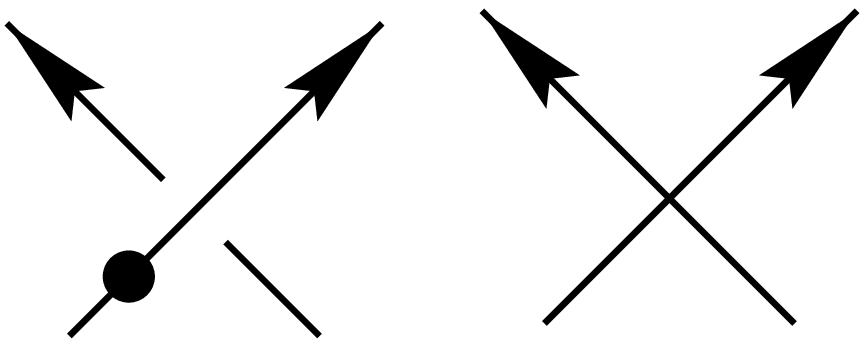,height=15pt}}) - \phi(\mbox{\epsfig{file=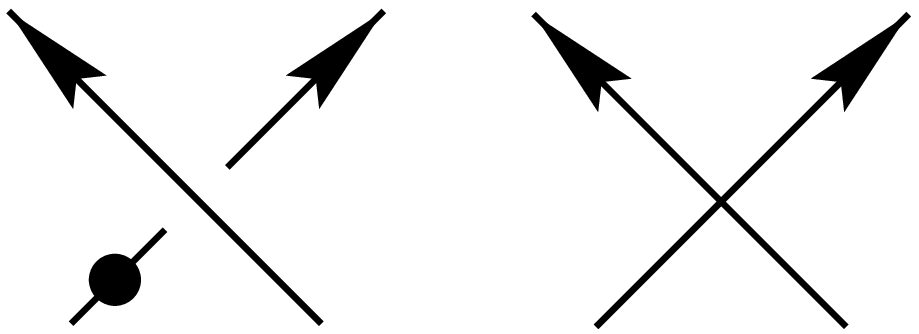,height=15pt}}) \\ \\
= \phi(\mbox{\epsfig{file=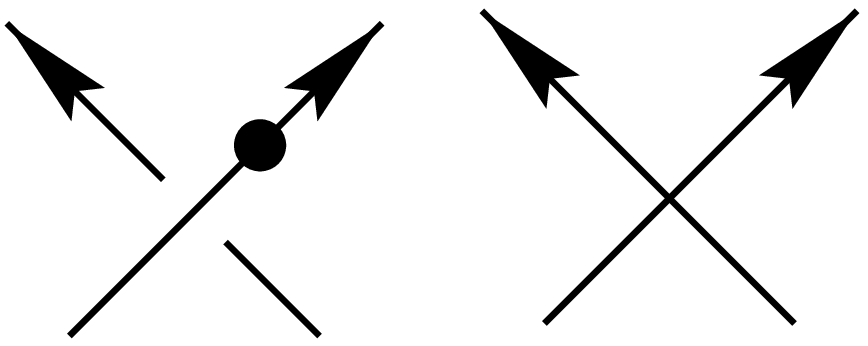,height=15pt}}) - \phi(\mbox{\epsfig{file=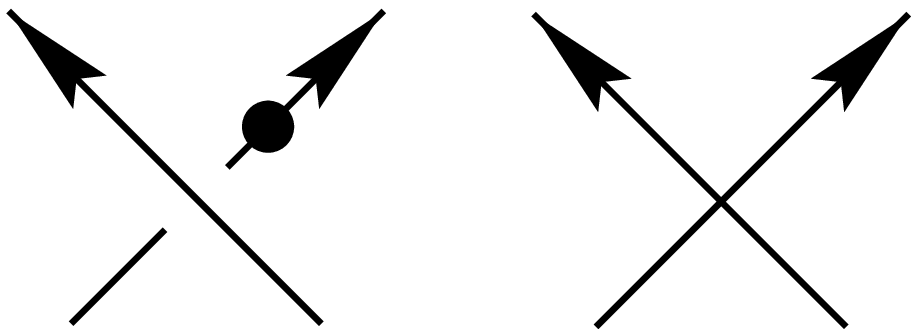,height=15pt}})
= \phi(\mbox{\epsfig{file=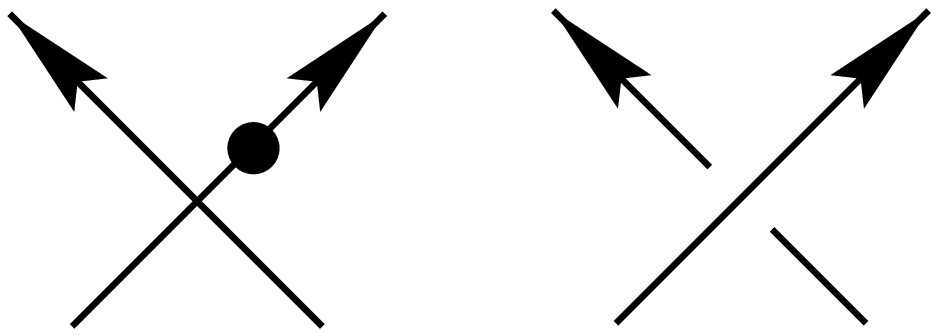,height=15pt}}) - \phi(\mbox{\epsfig{file=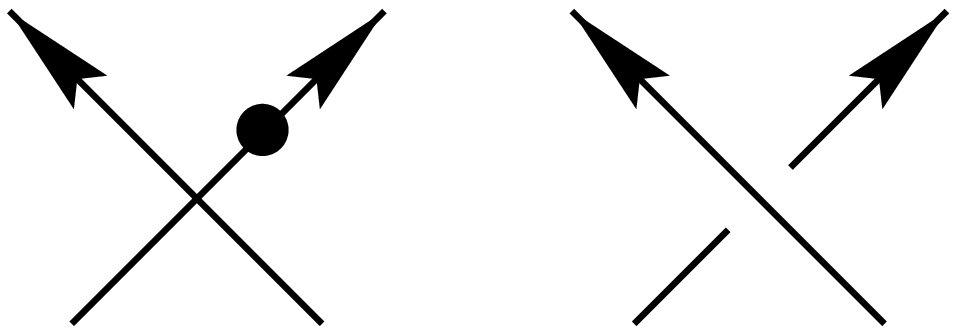,height=15pt}})
\end{array} \]
Define $\, \tilde{\phi}(\mbox{\epsfig{file=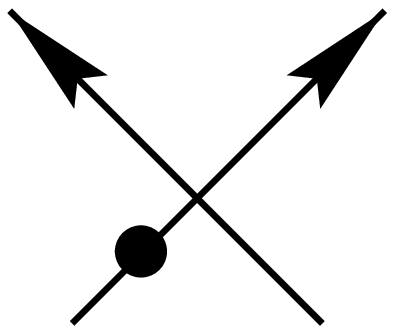,height=15pt}}) = \phi(\mbox{\epsfig{file=xb2.eps,height=15pt}}) - \phi(\mbox{\epsfig{file=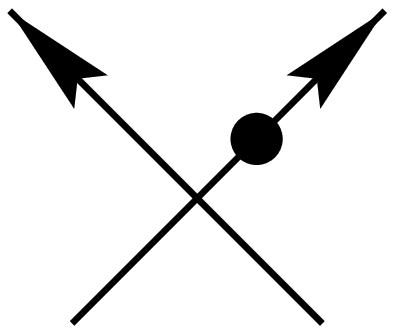,height=15pt}})$.
Clearly, $\phi$ is invariant under changes of base point if and only if $\tilde{\phi}=0$.
\[ \begin{array}{l}
\tilde{\phi}(\mbox{\epsfig{file=xp1.eps,height=15pt}})-\tilde{\phi}(\mbox{\epsfig{file=xm1.eps,height=15pt}}) = \\ \\
= \left(\phi(\mbox{\epsfig{file=xp1.eps,height=15pt}}) - \phi(\mbox{\epsfig{file=xp2.eps,height=15pt}})\right) - \left(\phi(\mbox{\epsfig{file=xm1.eps,height=15pt}}) - \phi(\mbox{\epsfig{file=xm2.eps,height=15pt}})\right) \\ \\
= \left(\phi(\mbox{\epsfig{file=xp1.eps,height=15pt}}) - \phi(\mbox{\epsfig{file=xm1.eps,height=15pt}})\right) - \left(\phi(\mbox{\epsfig{file=xp2.eps,height=15pt}}) - \phi(\mbox{\epsfig{file=xm2.eps,height=15pt}})\right) = 0
\end{array} \]
As $\tilde{\phi}$ is invariant under crossing changes and $\pi_1(\M)=1$, for all $K\in X_1^b,\; \tilde{\phi}(K)= \tilde{\phi}(\mbox{\epsfig{file=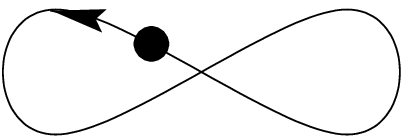,height=15pt}})= 0$. \eop


\begin{cor} \label{C2}
If $\pi_1(\M)=\pi_2(\M)=1,$ then the natural restriction map $\rho: H^1(CK)\to H^1(D)$ is an isomorphism.
\end{cor}
\pf\qquad
From the exact sequence of the pair $(CK,D)$ and Corollary \ref{C1}, it suffices to show that $\rho$ is surjective:
\[H^1(CK,D)=0\rightarrow H^1(CK)\overset{\rho}{\rightarrow} H^1(D) \rightarrow H^2(CK,D) \]
\[ \begin{array}{l}
Z^1(CK)=\{\phi\in C^1(CK): \delta\phi=0\}=\{\phi\in C^1(CK): \phi \text{ is differentiable} \} \\ \\
B^1(CK)=\{\phi\in C^1(CK): \phi=\delta\psi\}=\{\phi\in C^1(CK): \phi \text{ is integrable} \} \\
\end{array} \]
By Theorem \ref{integrable}, $H^1(CK)=\{\phi\in Z^1(CK): \phi$ does not satisfy the 1-term relation\}.

For any $K_{12}\in D_2, \; \del K_{12}=(-1)^{i+1}\di K_{12}$ since at least one crossing is degenerate.  If $\phi(\del D_2)=0$, then $\phi(\di K_{12})=0$; i.e., $\phi$ is invariant under crossing changes.  Since the chord diagram of degree 1 is unique, $\phi$ is constant.
\[ \begin{array}{l}
Z^1(D)=\{\phi\in C^1(D): \phi(\del D_2)=0\}=\{\phi\in C^1(D): \phi \text{ is constant} \} \\ \\
B^1(D)=\{\phi\in C^1(D): \phi=\delta\psi |_{D_1}=0 \} \\ \\
H^1(D)=Z^1(D) \\
\end{array} \]
For any $\phi\in Z^1(D)$, there exists an extension to the constant map $\hat{\phi}\in Z^1(CK)$ so that $\rho(\hat{\phi})=\phi$.
\eop
\begin{thm} {For $\M=\R^3, \; \V_2|_{C_2} \cong H^2(CK,D)$.} \label{H2eqV2}
\end{thm}
\pf\qquad From Corollary \ref{C2}, we obtain the following commutative diagram:
\[ \xymatrix{ {H^1(CK)} \ar[r]^\cong & {H^1(D)} \ar[r]^0 & {H^2(CK,D)} \ar[r] & {H^2(CK)} \ar[r] & {H^2(D)} \\
{} & {} & {} & {\V_2|_{C_2}} \ar[ul]^{\exists\gamma} \ar[u] \ar[ur]_0 & {} } \]
By Corollaries \ref{V2tilde} and \ref{H2bar} we obtain the following commutative diagram:
(Note by Corollary \ref{pi1}, $H^1(\tck)=0$.)
\[ \xymatrix{ {0} \ar[r] & {H^1(D\cup_fA)} \ar[r]^\alpha & {H^2(\tck,D\cup_fA)} \ar[r] & {H^2(\tck)} \ar[r]^0 & {H^2(D\cup_fA)} \\
{} & {} & {} & {\V_2|_{C_2}} \ar[ul]^{\exists\tilde{\gamma}} \ar[u]^\cong \ar[ur]_0 & {} } \]
We can relate these diagrams by the following commutative diagram, where the verticals are restriction maps, and the isomorphism is by excision:
\[ \xymatrix{ {H^1(D)} \ar[r]^0 & {H^2(CK,D)} \\
{H^1(D\cup_fA)} \ar[r]^\alpha \ar[u] & {H^2(\tck,D\cup_fA)} \ar[u]_\cong } \] \noindent Thus
$\alpha =0$, and consequently $\tilde{\gamma}:\V_2 \overset{\cong}{\to} H^2(\tck,D\cup_fA)
\cong H^2(CK,D)$. \eop

\subsection{Higher dimensional cohomology}

In higher dimensions, we consider the even and odd cases separately.  Our main problem is to find which Vassiliev invariants are nontrivial in cohomology, and a surprising disparity emerges in Corollary \ref{Cor2k}.  Moreover, in Corollary \ref{Conway} we show that Vassiliev invariants from the Conway polynomial--the best understood weight system of strictly even Vassiliev invariants--are all nontrivial in cohomology.

\begin{thm} \label{Vodd}
Suppose $\pi_1(\M)=\pi_2(\M)=1.$  For any Vassiliev invariant $\phi \in \V_N,$ \mbox{$[\phi]=0 \in H^{2k+1}(CK,D)$} $\quad \forall k \geq 0$.
\end{thm}

\pf\qquad
Of course, we must restrict $\V_N|_{C_{2k+1}}$.  Every Vassiliev invariant $\phi \in \V_N$ is integrable: \mbox{$\int(\phi|_{C_{2k+1}})= \phi|_{C_{2k}},$} so from the following lemmas we obtain that $[\phi]=0 \in H^{2k+1}(CK) \; \forall k \geq 0$.  Therefore, from (\ref{diag}), $\gamma(\phi) \in \ker \alpha$.  Since $\V_* \to H^*(D)$ is always zero, it follows that $[\phi]=0 \in H^{2k+1}(CK,D) \; \forall k \geq 0$.

\begin{lem} \label{lem1}
For $n=2k+1,$ if $\phi\in C^n(CK)$ is differentiable, then $\phi \in Z^n(CK).$
\end{lem}
\pf
\[ \delta\phi(K_{1 \ldots n+1})=\phi(\del K_{1 \ldots n+1})=\sum_{i=1}^{n+1} (-1)^{i+1} \phi(\di K_{1 \ldots n+1})=0 \]
because the alternating sum has an even number of equal terms.

\begin{lem} \label{lem2}
For $n=2k+1,$ if $\phi\in C^n(CK)$ is integrable, then $[\phi]=0 \in H^n(CK).$
\end{lem}
\pf\qquad
Since $\phi$ is integrable, it is also differentiable, so by the previous lemma, $\phi$ is an \n cocycle.  Let $\int\phi=\psi$.  Then
\[ \delta\psi(\Kn)=\sum_{i=1}^{n} (-1)^{i+1} \psi(\di \Kn)=\sum_{i=1}^{n} (-1)^{i+1} \phi(\Kn)=\phi(\Kn) \]
because the alternating sum has an odd number of equal terms, with the first and last terms positive.
\eop

Together with Theorem \ref{integrable}, if $\pi_1(\M)=\pi_2(\M)=1,$ these lemmas imply that for $n=2k+1,$ a differentiable invariant $\phi$ satisfies the 1-term and 4-term relations only if $\phi=\delta\psi$.  In this case, we obtain an explicit topological obstruction, that $\phi$ must be a coboundary, for the 1-term and 4-term relations to hold.

The even-dimensional case seems more interesting.  It is possible that all of the nontrivial cohomology of $CK$ arises from Vassiliev invariants of even order.

\begin{thm} \label{Veven}
Suppose $\pi_1(\M)=\pi_2(\M)=1.$  \; For any $n=2k, \; \V_n(\M) = \{\phi \in Z^n(CK,D) : \phi$ is integrable\}.
\end{thm}

\pf \qquad
Every Vassiliev invariant is an integrable cocycle by Proposition \ref{zn}.  Conversely, let $\phi$ be an integrable cocycle.  If $\phi$ is both integrable and invariant under crossing changes, then $\phi$ is a \emph{weight system}.  Since $\pi_1(\M)=1$, $\V_n(\M)$ contains a vector subspace isomorphic to $\V_n(\R^3)$ (Theorem 0.1 \c{Kalf}).  Over $\Q$, any given weight system can be integrated all the way to a knot invariant by the Kontsevich integral \c{Kont,BN}, so $\phi \in \V_n(\M)$.  The result now follows from the following lemma.

\begin{lem} \label{lem3}
For $n=2k$, suppose $\phi\in C^n(CK)$ is differentiable.  Then $\delta\phi=0$ if and only if $\phi$ is invariant under crossing changes.
\end{lem}

\pf \[ \delta\phi(K_{1 \ldots n+1})=\sum_{i=1}^{n+1} (-1)^{i+1} \phi(\di K_{1 \ldots n+1})=\phi(\d{n+1} K_{1 \ldots n+1}) \]
because the alternating sum has an odd number of equal terms, with the first and last terms positive.  Therefore, $\delta\phi=0$ if and only if $\phi(\Kp)-\phi(\Km)=\phi(\d{n+1}\Kx)=0.$
\eop

\begin{cor} \label{Cor2k}
If $\phi \in \V_N \setminus \V_{N-1}$ and $[\phi] \in H^n(CK,D)$ is nontrivial, then $N=n=2k$.
\end{cor}

\pf \qquad
By hypothesis, $\phi \in \V_N|_{C_n}$.  Since $[\phi] \neq 0, \; n \leq N$, and $n=2k$ by Theorem \ref{Vodd}.  By Theorem \ref{Veven}, $\phi \in \V_n|_{C_n}$, so $N \leq n$.  Therefore, $N=n=2k$.

 \eop

\begin{ques}
For $n=2k,$ suppose $\phi\in \V_n\setminus \V_{n-1}.$  Is $[\phi] \in H^n(CK,D)$ nontrivial?
\end{ques}

\begin{figure}
\begin{center}
\scalebox{.4}{\includegraphics{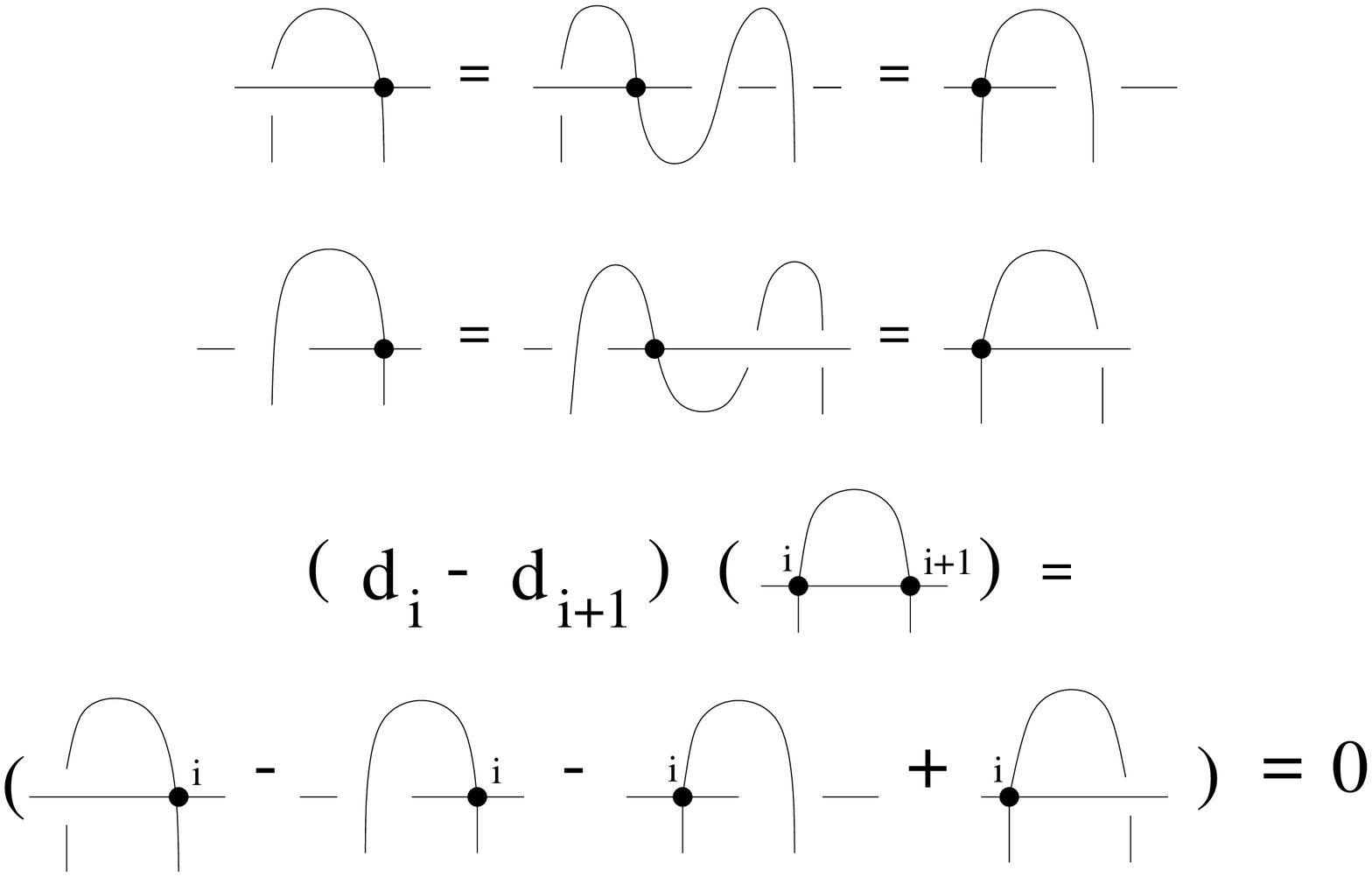}}
\caption{Proof of Theorem \ref{akira}} \label{akirafig}
\label{akira.fig}
\end{center}
\end{figure}

We have a partial result suggesting that the answer is yes, at least for $\R^3$.  The idea is the same as in Theorem \ref{V2}, to find for any chord diagram a representative $\Kn$, such that $\del \Kn=0$.  After considering this problem together, the following theorem was obtained by Akira Yasuhara.

\begin{thm}[Yasuhara] \label{akira}
For any $n \geq 1$, there exist Vassiliev invariants in $\V_{2n}(\R^3)$ which are nontrivial in $H^{2n}(CK(\R^3),D)$.
\end{thm}
\pf \qquad We can describe an ordered \n chord diagram as follows:  Let $C$ be a circle and let $p_1,p_2,\dots,p_{2n}$ points which lie on $C$ counterclockwise in this order. Let
\[ (i_1,j_1)(i_2,j_2)\ldots(i_n,j_n), \; i_k,j_k\in\{1,2,\dots,2n\} \]
denote an $n$-chord diagram obtained from $C$ by connecting $p_{i_k}$ and $p_{j_k}$, where $k=1,2,\dots,n$.  Now, let
\[ D_{2n}=(i_1,j_1)(i_2,j_2)\ldots(i_{2n},j_{2n}), \; \{i_1,i_2,\dots,i_{2n},j_1,j_2,\dots,j_{2n}\}=\{1,2,\dots,4n\} \]
be a $2n$-chord diagram, such that $|i_k-i_{k+1}|=|j_k-j_{k+1}|=1$ for any $k \in \{1,3,\dots,2n-1\}$.
In Figure \ref{akira.fig}, we show that there exists an ordered $2n$-singular knot $K_{1 \ldots 2n}$ which represents $D_{2n}$, such that the double point labeled by $k$ corresponds to the chord $(i_k,j_k)$, and $\left(\di-\d{i+1}\right)K_{1 \ldots 2n}=0$ for all $i\in \{1,3,\dots,2n-1\}$.  Therefore, $\del K_{1 \ldots 2n}=0$.

\begin{figure}
\begin{center}
\scalebox{.45}{\includegraphics{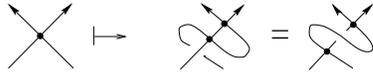}}
\caption{Perturbing an \n singular knot to a $2n$-singular knot}
\label{perturb}
\end{center}
\end{figure}

We can characterize such a $2n$-chord diagram as a perturbation of some \n chord diagram $D_n$ by replacing each chord of $D_n$ with an adjacent pair of chords, allowing either $||$ or $\times$ for each of these pairs.  In particular, any \n singular knot which represents $D_n$ can be perturbed to a $2n$-singular knot representing $D_{2n}$ by the local change shown in Figure \ref{perturb}.  Any $v \in \V_{2n}(\R^3)$ which is nonzero on linear combinations of these $2n$-singular knots is not a coboundary.  Therefore, $v$ represents a nontrivial cohomology class in $H^{2n}(CK(\R^3))$, and by (\ref{diag}) lifts to a nontrivial cohomology class in $H^{2n}(CK(\R^3),D)$.
\eop

\begin{cor} \label{Conway}
For any $v \in \V_{2n}(\R^3)$ coming from the Conway polynomial, $v$ is nontrivial in $H^{2n}(CK(\R^3),D)$.
\end{cor}
\pf \qquad Let $W_C$ be the Conway weight system.  For any \n chord diagram $D$, let $D'$ be the perturbed chord diagram obtained by replacing each chord of $D$ by an adjacent pair of intersecting chords.  By Theorem 2 of \c{BNG}, $W_C(D')=W_C(\bigcirc)=1$, where $\bigcirc$ denotes the chord diagram with zero chords.  Let $K_{1 \ldots 2n}$ represent $D'$ with any ordering such that every perturbed pair of double points is ordered consecutively modulo $2n$ (e.g., an ordering from any base point).  Now, $\del K_{1 \ldots 2n}=0$ by the proof of Theorem \ref{akira}, but $v(K_{1 \ldots 2n})=W_C(D')=1$, so $v$ is not a coboundary.
\eop

\begin{figure}
\begin{center}
\scalebox{.45}{\includegraphics{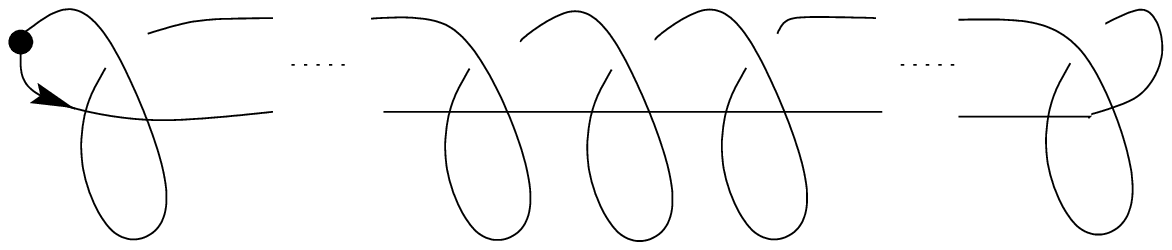}}
\caption{$2n$-singular knot $K_{1 \ldots 2n}$ with $\del K_{1 \ldots 2n}=0$ and $v_2^n(K_{1 \ldots 2n}) \neq 0$}
\label{v2.eps}
\end{center}
\end{figure}

\begin{ex} \label{V2n} \rm
For $v_2 \in \V_2(\R^3), \; v_2^n$ is nontrivial in $H^{2n}(CK(\R^3),D)$ because $v_2^n$ is nonzero on the $2n$-singular knot $K_{1 \ldots 2n}$ shown in Figure \ref{v2.eps}.  With notation as in Theorem \ref{akira}, $K_{1 \ldots 2n}$ represents the chord diagram $\{(i,4n-i)(i+1,4n-i+1): i=1,3,\dots,2n-1\}$.  With the base point as indicated, $\del K_{1 \ldots 2n}=0$.
\end{ex}

\begin{rmk} \rm
Any knot invariant $\phi \in C^0(CK)$ can be extended to an invariant of singular knots by the skein relation (\ref{skein}).  Then $\phi$ is integrable, but may not be of finite type.  Since the proof of Theorem \ref{Vodd} just requires integrability, we can show that if $\phi$ is not of finite type, then $[\phi]=0 \in H^n(CK,D), \; \forall n \geq 0$.
\end{rmk}
This leads naturally to the following question:

\begin{ques}
Is $H^{2k+1}(CK,D)$ nontrivial for any $k \geq 1$?
\end{ques}

\begin{rmk} \rm
Fenn, Rourke, and Sanderson \c{FRS} studied a similar construction with cubical sets and its relationship to Vassiliev invariants.  Indeed, our Lemmas \ref{lem1} and \ref{lem2} had appeared earlier in their paper (Lemma 7.2 \c{FRS}).  Their main result is that $\phi \in C^{n+1}(X)$ is integrable if and only if an obstruction in $H^1(J^n(X))$ vanishes, where $J^n(X)$ is the James complex of a $\square$-set $X$ (Theorem 7.4 \c{FRS}).  Their construction depends on ordering the double points with a base point, as in $CK_b$.  However, they did not consider other orderings, so that $CK_0$ may have different properties.  It would be interesting to relate our construction to their theory.
\end{rmk}

\subsection{Remarks on the mirror map} \label{mirror}

For $n\geq 0$ and any $\Kn$ in $X_n^b(\R^3)$ or $X_n^0(\R^3)$, let $\Kn^*$ be its mirror image.  For $X=X^b$ or $X^0$, respectively, let $f:C_n(X)\to C_n(X)$ be the map induced by the mirror map $\Kn \mapsto \Kn^*$.  We obtain:
\[ f\circ\di(\Kn)=(K_{\times\ldots+\ldots\times})^* - (K_{\times\ldots-\ldots\times})^*=-\di\Kn^*=-\di\circ f(\Kn) \]
Therefore, $f\del=-\del f$, so $f$ induces maps on both $H_*(CK)$ and $H^*(CK)$.  Clearly, if $\Kn$ is any degenerate knot, $f(\Kn)$ is also degenerate, so $f$ induces maps on the relative homology and cohomology.

Let $v_n$ be any canonical Vassiliev invariant (in the sense of \c{BNG}) of degree $n$.  By
the properties of the Kontsevich integral, $v_n(K^*)=(-1)^n v_n(K)$.  Consequently, for any
$v\in\V_n(\R^3)$, the Vassiliev invariant $\bar{v}(K)=v(K^*)-(-1)^n v(K)$ is of degree
strictly less than $n$.  By resolving double points, we obtain $v_n(\KM^*)=(-1)^{n+m}
v_n(\KM)$. By Corollary \ref{Cor2k}, if $v_n$ is nontrivial in $H^m(CK,D)$, then $n=m=2k$.
Thus, the induced map $f^*$ acts as the identity on the subgroup of $H^*(CK,D)$ generated by
canonical Vassiliev invariants, and hence on the subgroup generated by all Vassiliev
invariants.

On the other hand, we can define a topological map $F:CK(\R^3)\to CK(\R^3)$ obtained from the mirror map $\Kn \mapsto \Kn^*$.  As before, let $I^n=\{\vec{x} \in \R^n : -1 \leq x_i \leq 1 \}$.  For any $x\in I^n$, we define $F\circ\g(x) = g_{K_{1 \dots n}^*}(-x)$.  It is easy to verify that $F$ is well-defined.

Using the same notation for the induced map, $F: C_n(X)\to C_n(X)$ is given by $F(\Kn)=[F\circ\g(I^n)]$.  For $n\geq 0$, let $\r_n:I^n\to I^n$ be the map $\r_n(x)=-x$.  Since $\r_n$ is the product of reflections in each of the $n$ coordinates,
\[ F(\Kn)=[F\circ\g(I^n)] = [g_{K_{1 \dots n}^*}\circ\r_n(I^n)]= (-1)^n[g_{K_{1 \dots n}^*}(I^n)]= (-1)^n \Kn^* \]
As it is induced from a topological map, $F$ is a chain map:
\[ \begin{array}{l}
F(\del\Kn)=F\left(\sum\limits_{i=1}^n (-1)^{i+1} \di\Kn\right)=(-1)^{n-1}\sum\limits_{i=1}^n (-1)^{i+1}f\circ\di\Kn= \\ \\
= (-1)^{n-1}\sum\limits_{i=1}^n (-1)^{i+1} (-\di\Kn^*) = (-1)^n \del\Kn^*= \del F(\Kn)
\end{array} \]

As above, $F^*$ acts as the identity on the cohomology subgroup generated by Vassiliev
invariants. Moreover, $F$ acts by $(-1)^n$ on the cochain subgroup generated by canonical
Vassiliev invariants of degree $n$:
\[ v_n\circ F(\KM)=(-1)^m v_n(\KM^*)=(-1)^{n+2m}v_n(\KM)=(-1)^n v_n(\KM) \]

\begin{ques}  Do $f^*$ or $F^*$ act as the identity on $H^*(CK,D)$ or $H^*(CK)$?
\end{ques}



Geometrically, we can say that $F$ is a bijection on \n cells of $CK$.   Gillete and Van Buskirk \c{GVB} described several related examples, including a knot diagram of minimal crossing number with a crossing which can be switched to obtain the mirror image knot.  This provides an instance when $F$ maps an edge to itself, but with reversed orientation.

\begin{conj}{If $\Phi:CK(\R^3) \overset{\cong}{\to} CK(\R^3)$ is a bijection on vertices
which preserves edges then $\Phi$ is the identity map or the map $F$.}
\end{conj}

\begin{rmk} \rm
We say $\phi \in C^n(CK)$ is \emph{even} if $\phi \circ f=\phi$, and $\phi$ is \emph{odd} if $\phi \circ f=-\phi$.  Willerton \c{W2} gave an explicit formula showing that any even integrable invariant has a unique odd integral $\psi=\int\phi$, and that any integral of an odd invariant is even.
\end{rmk}

\subsection{Cup product in $H^*(CK)$} \label{cup}

Cup products for cubical cohomology have been discussed in standard texts (see e.g., \c{HW}).  Let $Q_n(X)$ be the free abelian group generated by singular \n cubes $g: I^n\to X$, and let $D_n(X)$ be generated by degenerate \n cubes.  The main difficulty is finding a chain map $\Delta: C_n(X)\to(C_*(X)\otimes C_*(X))_n$, where $C_n(X)=Q_n(X)/D_n(X)$.  In Section 9.3 of \c{HW}, the map $\Delta$ given below is shown to satisfy all of the required properties.  Let $H=(h_1\ldots h_p)\subset\{1\ldots n\}$ be any subset (possibly empty) with the natural order on its elements, and let $K$ be the complementary subset with the natural order on its elements.  Let $\rho_{(HK)}$ denote the sign of $\s\in S_n$, where $\s(HK)=(1\ldots n)$.  For $\e\in\{\pm 1\}$, let $\lambda^\e_H(u_1,\dots,u_p)=(v_1,\dots,v_n)$, where $v_i=\e$ if $i\notin H$ and $v_{h_r}=u_r, \, r=1,\dots,p.$  Thus, $\lambda^{-1}_H$ is an isometry of $I^p$ onto a particular back $p$-face of $I^n$, and $\lambda^{+1}_H$ maps onto the parallel front $p$-face of $I^n$.  In our context, $C_n(CK)$ is generated by maps $\g:I^n\to CK$.  We define $\Delta: C_n(CK) \to (C_*(CK)\otimes C_*(CK))_n$ by

\[ \begin{array}{l}
\Delta(\g)=\sum\limits_{H\subset \{1\ldots n\}} \rho_{(HK)} \g\circ\lambda^{-1}_H \otimes \g\circ\lambda^{+1}_K \\ \\
=\sum\limits_H \rho_{(HK)} g_{K_{\times\ldots a_1\ldots a_q\ldots\times}}\otimes g_{K_{\times\ldots b_1\ldots b_p\ldots\times}},\text{ where } a_i=-1\text{ and } b_i=+1 \; \forall i \\ \\
=\sum\limits_{H\subset \{1\ldots n\}} \rho_{(HK)} g_{\d{K}^- \Kn} \otimes g_{\d{H}^+ \Kn} \\
\end{array} \]
In our notation, $\d{H}^+$ takes the positive resolution of double points with labels in $H$, and similarly for $\d{K}^-$.  We can simply refer to the singular knots:
\begin{equation}
\Delta(\Kn)=\sum_{H\subset \{1\ldots n\}} \rho_{(HK)}\d{K}^- \Kn\otimes\d{H}^+ \Kn
\end{equation}
For any $u\in H^p(CK)$ and $v\in H^q(CK), \; u\cup v=\Delta^*(u\times v)\in H^{p+q}(CK)$, and similarly for relative cohomology.  Therefore,
\[ u\cup v(\Kn) = \sum_{H\subset \{1\ldots n\}} \rho_{(HK)}u(\d{K}^-\Kn)\cdot v(\d{H}^+\Kn) \]
Bar-Natan has suggested that this expression is similar to one arising in a variant of graph cohomology related to Vassiliev invariants.

\begin{ex}\rm
By Theorem \ref{V2}, $v_2\in\V_2$ is a generator of $H^2(CK(\R^3))$.  To compute $v_2\cup v_2(K_{1234})$, note that $v_2(\d{K}^- K_{1234})\cdot v_2(\d{H}^+ K_{1234}) \neq 0$ only for $p=2$.
\[ \begin{array}{ccccc}
H & K & \rho_{(HK)} & \d{K}^- K_{1234} & \d{H}^+ K_{1234} \\ \\
1,2 & 3,4 & + & K_{\times\times - -} & K_{+ + \times\times} \\
1,3 & 2,4 & - & K_{\times - \times -} & K_{+\times +\times} \\
1,4 & 2,3 & + & K_{\times - - \times} & K_{+ \times\times +} \\
2,3 & 1,4 & + & K_{- \times\times -} & K_{\times + + \times} \\
2,4 & 1,3 & - & K_{- \times - \times} & K_{\times + \times +} \\
3,4 & 1,2 & + & K_{- - \times\times} & K_{\times\times + +} \\
\end{array} \]
\[ \begin{array}{l}
v_2\cup v_2(K_{1234}) = \\ \\
v_2(K_{\times\times - -})v_2(K_{+ + \times\times})-v_2(K_{\times - \times -})v_2(K_{+\times +\times})+v_2(K_{\times - - \times})v_2(K_{+ \times\times +}) \\ \\
+v_2(K_{- \times\times -})v_2(K_{\times + + \times})-v_2(K_{- \times - \times})v_2(K_{\times + \times +})+v_2(K_{- - \times\times})v_2(K_{\times\times + +}) \\
\end{array} \]
As we discuss in Section \ref{Hopf}, $v_2\cdot v_2(K_{1234})$ is equal to the same expression but without the signs $\rho_{(HK)}$ \c{W1}.
\end{ex}

\begin{rmk} \rm
On page 210 of \c{Q}, Quillen shows that if $A=\oplus_{n\geq 0}A_n$ is a graded commutative
algebra over $\Q$ with $A_n$ finite dimensional for all $n$, and $A_0=\Q, A_1=0$, then $A$ is
the rational cohomology ring of a simply-connected pointed topological space.  The algebra of
Vassiliev invariants under the cup product satisfies these conditions.  Let $\V=(\oplus_{n\geq
0}\V_n,\cup)$.  Therefore, there exists a simply-connected pointed space $X$ such that
$H^*(X,\Q)=\V$.  For even Vassiliev invariants, the simply connected knot complex $\tck(\R^3)$ seems very close to such a space
since by Corollary \ref{H2bar} and Theorem \ref{H2eqV2}, $\V_2\cong H^2(\tck)\cong H^2(CK,D)$,
and by (\ref{H2surj}), $H^p(\tck)\cong H^p(CK) \; \forall p\geq 4$.
\end{rmk}

\section{Hopf algebras of ordered chord diagrams} \label{Hopf}

In this section, we show how the cup product in $H^*(CK_0)$ for Vassiliev invariants arises naturally from a bialgebra structure on \emph{ordered} chord diagrams.  The resulting Hopf algebra $\aa^0$ has an additional differential structure:  $\aa^0$ is a cocommutative differential graded Hopf algebra.  We also consider the quotient algebra of \emph{oriented} chord diagrams, $\aa^\omega$, a differential graded Hopf algebra which is commutative and cocommutative in the graded sense.

\subsection{Additional structures on chord diagrams}

We follow the notation of Section \ref{structures}.  Let $\dd$ denote the free abelian group over $\Q$ generated by chord diagrams with each boundary circle directed counterclockwise.  The \emph{degree} of a chord diagram is the number of chords.  Let $\dd_n$ denote the subgroup generated by chord diagrams of degree $n$.  We say $D$ is an \emph{ordered chord diagram} if its chords are ordered from $1$ to $n$.  Two ordered chord diagrams are equivalent if one can be obtained from the other by some rotation.  Let $\dd^0$ denote the free abelian group over $\Q$ generated by equivalence classes of ordered chord diagrams.

Let $\dd^b$ denote the free abelian group over $\Q$ generated by chord diagrams with a base point.  $\dd^b$ is equivalent to the group generated by \emph{linearized chord diagrams}, where chords are identifications of points on the directed real line.  For any chord diagram in $\dd^b$, the chords can be ordered by their left endpoints on the real line.  Let $\dd_n^{b,0}$ be generated by elements in $\dd_n^b$ with the ordering of the chords permuted by any $\s\in S_n$.  Let $\dd_0^{b,0}=\dd^b_0$.  As sets $\dd_n^{b,0} \cong \dd^b_n \times S_n$.

Just as we defined an orientation on singular knots, we can consider \emph{oriented chord diagrams}.  For any $\Dn\in \dd_n^0$ and any $\sigma\in S_n$, let $\Dsn$ denote the chord diagram obtained by permuting the ordering of the chords of $\Dn$.  Let $\dd_n^\omega$ denote the abelian group generated by elements $\Dn\in\dd_n^0,$ subject to the following relations:
\begin{equation} \label{Dw}
\Dn \sim_\omega \text{\rm sign} (\s) \Dsn \; \text{ for any } \s \in S_n
\end{equation}

As in the diagram (\ref{CnX}), by taking the corresponding projections the following diagram commutes:
\begin{equation} \label{dd}
\xymatrix{ {} & \ar[dl]_{p_b} {\dd^{b,0}} \ar[dr]^{p_0} & {}  \\
{\dd^b} \ar[dr]_{p_\omega\circ i}  & {} & \ar[dl]^{p_\omega} {\dd^0} \\
{} & {\dd^\omega} & {} }
\end{equation}

\subsection{Hopf algebra of chord diagrams }

Let $\aa$ be the usual Hopf algebra of chord diagrams modulo the 4T relation.  We denote its bialgebra structure by $\aa=(\aa, \cdot, \D)$.  The product $D^1\cdot D^2$ is obtained from the direct sum operation on knots, and is well-defined as a consequence of the 4T relation.  This product is commutative and associative.  The coproduct, which is cocommutative and coassociative, is obtained from the usual shuffle coproduct on a tensor algebra (dual to the tensor product):

\begin{defn} \label{Delta}
We define $\D : \aa \to \aa \otimes \aa$ as follows:  For $D\in\dd_n$, we choose any ordering of its chords from $1$ to $n$.  Let $H=(h_1\ldots h_p)\subset\{1\ldots n\}$ be any subset (possibly empty) with the natural order on its elements.  Let $K$ be the complementary subset with the natural order on its elements.  Let $D_H$ denote the chord diagram obtained from $D$ by removing chords with labels in $H$.  Then
\[ \D(D) =  \sum_{H\subset \{1\ldots n\}} D_K \otimes D_H \]
\end{defn}

If $v\in \V_m$ and $D\in \dd_m$, then $W_m(v)(D)=v(K_D)$ defines a weight system.  The following algebraic relations are well known (see, e.g., \c{BN, W1}).  Let $v_1\in \V_p$ and $v_2\in \V_q$. Let $n=p+q$.
\[ \begin{array}{c}
\D(D^1\cdot D^2)=\D(D^1)\D(D^2) \\ \\
v_1\cdot v_2(\Kn) = \sum\limits_{H\subset \{1\ldots n\}} v_1(\d{K}^-\Kn)\cdot v_2(\d{H}^+\Kn) \\ \\
W_n(v_1\cdot v_2) = (W_p(v_1)\otimes W_q(v_2))\circ \D
\end{array} \]

For linearized chord diagrams, there is a corresponding linearized 4T relation (also known as the STU relation).  Modulo this relation, any mutually disjoint sets of chords can be moved past each other.  Therefore, the product $D^1\cdot D^2$ given by concatenation is preserved by the natural map $\dd_n^b\to\dd_n$, which is induced by the map $X_n^b\to X_n$.  Let $\aa^\ell$ be the quotient of $\dd^b$ by the linearized 4T relation.

\begin{prop} [Theorem 6, Lemma 3.1 \c{BN}] \label{Al}
$\aa^\ell \cong \aa$.
\end{prop}

\subsection{Hopf algebra of ordered chord diagrams } \label{Hopford}

We define a bialgebra structure on \emph{ordered} chord diagrams which is compatible with the cup product for Vassiliev invariants:
\[ v_1\cup v_2(\Kn) = \sum_{H\subset \{1\ldots n\}} \rho_{(HK)}v_1(\d{K}^-\Kn)\cdot v_2(\d{H}^+\Kn) \]

\begin{figure} 
\begin{center}
\scalebox{.55}{\includegraphics{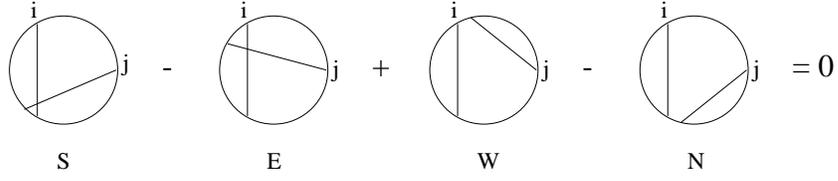}}
\caption{Ordered 4T relation}
\label{ord4T}
\end{center}
\end{figure}

The 4T relation is given by four diagrams, which are the same except for one ``fixed'' chord and one ``moving'' chord.  The \emph{ordered 4T relation} on $\dd^0$ is given by the same expression, where the fixed chord and the moving chord have the same label in all four diagrams.  (See Figure \ref{ord4T}.)  Thus, for each 4T relation on $\dd_n, \, n\geq 2$, we obtain $n(n-1)$ ordered 4T relations on $\dd_n^0$.

Let $\aa^0$ be the quotient of $\dd^0$ by all ordered 4T relations.  We will denote its bialgebra structure by $\aa^0=(\aa^0, \cup, \D^0)$.  If $D^1 \in \dd_p^0$ and $D^2 \in \dd_q^0$, the product $D^1\cup D^2$ is defined to be the chord diagram $D^1\cdot D^2$ with its ordering given by the same labels for the chords from $D^1$ and by labeling the $i^{\text{th}}$ chord of $D^2$ by $p+i$.

\begin{prop}
The cup product on $\aa^0$ is well defined.
\end{prop}
\pf \qquad The same argument used to prove Proposition 4.4 in \c{W1} will show the cup product is well defined on $\dd^0$ modulo the ordered 4T relation.   \eop

\noindent
Therefore, $\aa^0 \otimes \aa^0$ is an algebra with the following product:
\[ (a_1\otimes a_2)\cup (b_1\otimes b_2) = (-1)^{|a_2||b_1|}(a_1\cup b_1)\otimes (a_2\cup b_2) \]

For any $D \in \dd^0$, let $D_H$ be obtained by removing the chords with labels in $H$, with its ordering induced from $D$.  Let $\rho_{(HK)}$ be defined as in our cup product formula for $H^*(CK)$.
\begin{defn}
We define $\D^0 : \dd_n^0 \to \dd_n^0 \otimes \dd_n^0$ as follows:
\[ \D^0(D) =  \sum_{H\subset \{1\ldots n\}} \rho_{(HK)} D_K \otimes D_H \]
\end{defn}
It is easy to see that $\D^0$ descends to a coproduct on $\aa_n^0$:
\[ \D^0(\text{\rm ordered 4T})= (\text{\rm ordered 4T})\otimes\bigcirc + \bigcirc\otimes(\text{\rm ordered 4T}) \]
where $\bigcirc$ denotes any chord diagram in the ordered 4T relation with the chords $i$ and $j$ removed.

Let $T:\aa^0 \otimes \aa^0 \overset{\cong}{\rightarrow} \aa^0 \otimes \aa^0$ be the graded flip map, given by
\[T(a\otimes b)=(-1)^{|a||b|}b\otimes a \]
We can easily see that $\D^0$ is \emph{cocommutative} (i.e., $T\circ\D^0=\D^0$):
\[ T\circ\D^0(D) = \sum_{|H|=p} (-1)^{pq}\rho_{(HK)} D_H \otimes D_K=\sum_H \rho_{(KH)} D_H \otimes D_K=\D^0(D) \]

\begin{prop} \label{Dcup}
For $D^1, D^2 \in \aa^0, \quad \D^0(D^1\cup D^2)=\D^0(D^1)\cup\D^0(D^2)$.
\end{prop}
\pf \qquad Suppose $D^1 \in \aa_p^0$ and $D^2 \in \aa_q^0$.  Let $n=p+q$.
\[ \begin{array}{l}
\D^0(D^1)\cup\D^0(D^2)= \\ \\
= \left(\sum\limits_{H_1\subset \{1\ldots p\}} \rho_{(H_1K_1)} D^1_{K_1} \otimes D^1_{H_1}\right)\cup\left(\sum\limits_{H_2\subset \{1\ldots q\}} \rho_{(H_2K_2)} D^2_{K_2} \otimes D^2_{H_2}\right) \\ \\
= \sum\limits_{H_1, H_2}\rho_{(H_1K_1)}\rho_{(H_2K_2)}(-1)^{|D^1_{H_1}||D^2_{K_2}|}D^1_{K_1}\cup D^2_{K_2}\otimes D^1_{H_1}\cup D^2_{H_2} \\ \\
= \sum\limits_{H\subset \{1\ldots n\}}\rho_{(HK)}(D^1\cup D^2)_K\otimes (D^1\cup D^2)_H = \D^0(D^1\cup D^2)
\end{array} \]
The third equality follows from the fact that if $H=(H_1,H_2)$ and $K=(K_1,K_2)$, then $\rho_{(H_1K_1)}\rho_{(H_2K_2)}(-1)^{|K_1||H_2|}=\rho_{(HK)}$. \eop

If $v\in \V_n$ and $\Dn \in \aa_n^0$ then $W_n(v)(\Dn)=v(\Kn)$ defines a weight system in $(\aa^0_n)^*$, where $W_n:\V_n\to(\aa^0_n)^*$ is a graded map.  Thus by convention,
\[(W_p\otimes W_q)(v_1\otimes v_2)=(-1)^{pq}W_p(v_1)\otimes W_q(v_2)\]

\begin{thm} \label{Wcup}
Let $\mu(v_1\otimes v_2)=v_1\cup v_2$.  Then,
\[W_{p+q}\circ\mu = (\D^0)^*\circ (W_p\otimes W_q) \]
i.e., if $v_1\in \V_p$ and $v_2\in \V_q$ then \[ W_{p+q}(v_1\cup v_2) = (-1)^{pq}(W_p(v_1)\otimes W_q(v_2))\circ \D^0 \]
\end{thm}
\pf \qquad Let $n=p+q$.  For $\Dn\in\dd_n^0$, let $\Kn$ be any representative \n singular knot.
\[ \begin{array}{lr}
W_n(v_1\cup v_2)(\Dn) = v_1\cup v_2(\Kn) & \\ \\
= \sum\limits_{H\subset \{1\ldots n\}} \rho_{(HK)}v_1(\d{K}^-\Kn)\cdot v_2(\d{H}^+\Kn) & \\ \\
= \sum\limits_{|H|=p} \rho_{(HK)}v_1(\d{K}^-\Kn)\cdot v_2(\d{H}^+\Kn) & \\ \\
= \sum\limits_{|H|=p} \rho_{(HK)}W_p(v_1)(D_K)\cdot W_q(v_2)(D_H) & \\ \\
= \sum\limits_{H\subset \{1\ldots n\}} \rho_{(HK)}W_p(v_1)(D_K)\cdot W_q(v_2)(D_H) & \\ \\
= (-1)^{pq}(W_p(v_1)\otimes W_q(v_2))\circ \D^0(\Dn) & \hspace*{1cm} \square
\end{array}  \]

\begin{rmk} \rm
Recall that chord diagrams with a base point are equivalent to linearized chord diagrams.  Previously, we defined $\aa^\ell$ to be the quotient of $\dd^b$ by the linearized 4T relation.  Let $\aa^{b,0}$ be the quotient of $\dd^{b,0}$ by all linearized ordered 4T relations.  We can easily extend the results of this subsection to show a bialgebra structure $(\aa^{b,0}, \cup, \D^{b,0})$.  Bar-Natan pointed out that the proof of Proposition \ref{Al} can be modified by using the ordered 4T relations to prove $\aa^{b,0} \cong \aa^0$.
\end{rmk}

\subsection{Hopf algebra of oriented chord diagrams }

Let $\aa^\omega$ be the quotient $\aa^0/\sim_\omega$, with the relation $\sim_\omega$ as in (\ref{Dw}).
By the following proposition, the cup product on $\aa^0$ descends to a product on $\aa^\omega$.  Thus, we will denote its bialgebra structure by $\aa^\omega=(\aa^\omega, \cup, \D^\omega)$.

\begin{prop}
The cup product on $\aa^\omega$ is well defined.
\end{prop}
\pf \qquad Suppose $D^1 \in \aa_p^0$ and $D^2 \in \aa_q^0$.  If $\s\in A_p$ and $\t\in A_q$, then these act on disjoint index sets of $D^1\cup D^2$, and $\s\cdot\t\in A_{p+q}$.  Therefore,
\[ D^1_{\s(1\dots p)}\cup D^2_{\t(1\dots q)} \sim_\omega D^1_{1\dots p}\cup D^2_{1\dots q} \]
 \eop
The cup product on $\aa^\omega$ is commutative in the graded sense:
\[ D^1\cup D^2 = (-1)^{pq} D^2\cup D^1 \]
Therefore, the cup product on $\aa^\omega \otimes \aa^\omega$ is also graded commutative.

By the following proposition, the coproduct $\D^0$ on $\aa^0$ descends to a coproduct $\D^\omega$ on $\aa^\omega$.
\begin{prop} \label{Deltaw}
For $\Dn\in\dd_n^0,\,\s\in S_n, \; \D^0(\Dsn) \sim_\omega$ {\rm sign}$(\s)\D^0(\Dn)$.
\end{prop}
\pf \qquad For any $H\subset \{1\ldots n\}$, suppose $H=\{h_1,\dots,h_p\}$, where $h_i<h_{i+1}$ for $1\leq i\leq p-1$.  Let $\s(H)$ denote the ordered set $(\s(h_1),\dots,\s(h_p))$.  Let $\s[H]$ denote the ordered set $(\s(h_{i_1}),\dots,\s(h_{i_p}))$, where $\s(h_{i_j})<\s(h_{i_{j+1}})$ for $1\leq j\leq p-1$.  Let $\t_H\in S_p$ be the permutation which takes $\s[H]$ to $\s(H)$.  For $|H|\leq 1$, let sign$(\t_H)=+1$.
We obtain
\[ \rho_{(\s[H]\s[K])}\cdot\text{\rm sign}(\t_H)\cdot\text{\rm sign}(\t_K) =  \rho_{(\s(H)\s(K))} = \text{\rm sign}(\s)\cdot\rho_{(HK)} \]
Therefore,
\[ \begin{array}{l}
\D^0(\Dsn) = \sum\limits_{H\subset \{1\ldots n\}} \rho_{(HK)} (\Dsn)_K \otimes (\Dsn)_H \\ \\
= \sum\limits_{H\subset \{1\ldots n\}} \rho_{(\s[H]\s[K])} (\Dsn)_{\s[K]} \otimes (\Dsn)_{\s[H]}\\ \\
\sim_\omega \sum\limits_{H\subset \{1\ldots n\}}\rho_{(\s[H]\s[K])}\cdot\text{\rm sign}(\t_H)D_K \otimes \text{\rm sign}(\t_K) D_H \\ \\
= \sum\limits_{H\subset \{1\ldots n\}} \text{\rm sign}(\s) \rho_{(HK)} D_K \otimes D_H = \text{\rm sign}(\s)\D^0(\Dn)
\end{array} \]
 \eop

\begin{cor}
For $D^1, D^2 \in \aa^\omega, \quad \D^\omega(D^1\cup D^2)=\D^\omega(D^1)\cup\D^\omega(D^2)$.
\end{cor}
\pf \qquad  This follows from the Propositions \ref{Dcup} and \ref{Deltaw}.
\eop

By the same argument as for $\D^0$ above, $\D^\omega$ is graded cocommutative, so $\aa^\omega$ is a graded commutative and cocommutative Hopf algebra.

Because of the symmetry of chord diagrams, modulo the ordered 4T relations, often $\Dn=\Dsn$ for an odd permutation $\s$; if so, $\Dn \sim_\omega 0$.  All chord diagrams in degrees two and three are killed in this way.

\begin{ques} \label{aomega} Is $\aa^\omega$ nontrivial?
\end{ques}

\subsection{Differential structure on chord diagrams }

We now define the differential structure on $\aa^0$ and $\aa^\omega$.  For any $\Dn \in \dd^0$, we define
\[ \d{i}\Dn = \Dn \setminus \{i^{\text{th}}\;\text{chord}\} \]
with its ordering induced from $\Dn$.

\begin{defn}
For any $\Dn \in \dd^0$, define the boundary operator by
\[ \del \Dn = \sum_{i=1}^{n} (-1)^{i+1} \d{i} \Dn \]
\end{defn}

By the same arguments as in Propositions \ref{del2} and \ref{del3}, we obtain that $\del^2 = 0$ and $\del D_{\sigma(1 \cdots n)}=\text{sign} (\sigma)\del\Dn$.  Therefore, by the following proposition $\del: \aa_n^0 \to \aa_{n-1}^0$ and $\del: \aa_n^\omega \to \aa_{n-1}^\omega$ are well-defined.

\begin{prop}{$\del(${\rm ordered 4T}$)=0.$}
\end{prop}
\pf \qquad  For any chord diagrams $S,E,W,N$ as in Figure \ref{ord4T}, which are the same except for the part of the diagram shown, we wish to show $\del(S-E+W-N)=0$.
For $k \neq i,j, \; \d{k}(S-E+W-N)=0$ by the ordered 4T relation, so we obtain the following expressions:
\[ \del S= (-1)^i \d{i}S + (-1)^j \d{j}S  = (-1)^i \d{i}N + (-1)^j \d{j}N  = \del N \]
\[ \del E= (-1)^i \d{i}E + (-1)^j \d{j}E  = (-1)^i \d{i}W + (-1)^j \d{j}W  = \del W \]
 \eop
\begin{prop} $\del$ is a derivation with respect to the cup product.
\end{prop}
\pf \qquad Let $D^1, D^2 \in \dd^0$ be diagrams of degree $p$ and $q$, respectively.
\[ \begin{array}{l}
\del (D^1\cup D^2)= \sum\limits_{i=1}^{p+q} (-1)^{i+1} \d{i}(D^1\cup D^2) \\ \\
= \sum\limits_{i=1}^{p} (-1)^{i+1} \d{i}(D^1\cup D^2) + \sum\limits_{i=p+1}^{q} (-1)^{i+1} \d{i}(D^1\cup D^2) \\ \\
= (\del D^1)\cup D^2 + (-1)^p D^1\cup(\del D^2)
\end{array} \]
\quad \eop

The differential structure on $\aa^0 \otimes \aa^0$ is given by $\DD=\del\otimes 1 + 1\otimes\del$, so that $\DD(a \otimes b)=\del a\otimes b + (-1)^{|a|} a \otimes \del b$.  Let $\mu(a \otimes b)=a\cup b$.  We have proved
\begin{equation} \label{deriv}
\del\circ\mu=\mu\circ\DD
\end{equation}

Let $\aa^*$ be the graded dual of $\aa^0$, so that $\aa^*=\oplus \aa_n^*$.  Let $\delta=\del^*$.  The differential structure on $\aa^* \otimes \aa^*$ is given by $\DD^*=\delta\otimes 1 + 1\otimes\delta$.  By dualizing (\ref{deriv}), we obtain that $\delta$ is a \emph{coderivation}; i.e.,
\[ \D^0\circ\delta=\DD^*\circ\D^0 \]
and similarly for $\D^\omega$.

\begin{rmk} \rm
In Theorem \ref{akira}, we showed that for any $D\in\dd_n$, if $D'$ is obtained by any perturbation, then there exists a knot $K_{D'}\in X_{2n}^0$ which represents $D'$, such that $\del K_{D'}=0$.  Now, we can see that with the same ordering, $\del D'=0$.
\end{rmk}

